\documentclass[11pt] {article}      
\usepackage{amsmath,amssymb, amsfonts, amsthm, eucal}  
\usepackage{latexsym, graphicx, amscd}

\usepackage{eso-pic,graphicx,transparent}

\usepackage{newtxtext,newtxmath} 

\usepackage{enumitem}

\usepackage[T1]{fontenc}     
\usepackage{enumerate}       
\usepackage{float}                 
\usepackage{tikz}
\usetikzlibrary{matrix, arrows, calc}

\usepackage{mathtools}  
\usepackage{relsize}      %

\usepackage{tikz}
\usetikzlibrary{arrows,shapes}
\usetikzlibrary{chains}
\usetikzlibrary{matrix}
\definecolor{gold}{rgb}{1,.70,.0}   
\usepackage{float}                 

\theoremstyle{plain}
\newtheorem{theorem}{Theorem}[section]            

\theoremstyle{definition}

\numberwithin{theorem}{section}
\numberwithin{equation}{section}
\newcommand{\dotfraction}[2]{\genfrac{}{}{0.5pt}{}{#1}{#2}%
                        \!\lower.5pt\hbox{{$\circ$}}}

\def\Prm{\rm P}

\usepackage{titlesec}                                                           
\titlelabel{\thetitle.\ }                                                         
\titleformat*{\section}{\fontsize{14pt}{14pt} \bf}        
\usepackage{fancyhdr}
\usepackage{sidecap}  

\usepackage[hang,small,sc]{caption}

\let\OLDthebibliography\thebibliography
\renewcommand\thebibliography[1]{
  \OLDthebibliography{#1}
  \setlength{\parskip}{3pt}
  \setlength{\itemsep}{0pt plus 0.3ex}
}

\newcommand{\QED}{\ $\square$}

\def\vvec#1#2{
{\renewcommand*{\arraystretch}{.7}
\begin{bmatrix}#1\\#2\end{bmatrix}}
}

\makeatletter
\newcommand*\bigcdot{\mathpalette\bigcdot@{.41}}
\newcommand*\bigcdot@[2]{\mathbin{\vcenter{\hbox{\scalebox{#2}{$\m@th#1\bullet$}}}}}
\makeatother
\makeatletter
\newcommand*\bbigcdot{\mathpalette\bigcdot@{.61}}
\newcommand*\bbigcdot@[2]{\mathbin{\vcenter{\hbox{\scalebox{#2}{$\m@th#1\bullet$}}}}}
\makeatother

\def\curlrm{\hbox{\rm curl}\,}
\def\divrm{\hbox{\rm div}\,}

\def\SL{\hbox{\rm SL}}
\def\SO{\hbox{\rm SO}}

\def\spin{\hbox{\sf spin}}

\def\symp{\bigstar}  
\def\symp{\scalebox{.77}{$\star$}}  

\newcommand{\gaction}[2]{\genfrac{}{}{0.5pt}{}{#1}{#2}%
                       \!\!\lower2pt\hbox{\rotatebox[origin=c]{-90}{{$\looparrowright$}}}}
\newcommand{\ggaction}[2]{%
                        \lower2pt\hbox{\rotatebox[origin=c]{90}{{$\looparrowleft$}}}\!\!
                        \genfrac{}{}{0.5pt}{}{#1}{#2}%
                       \!\!\lower2pt\hbox{\rotatebox[origin=c]{-90}{{$\looparrowright$}}}
}

\newcommand{\ogaction}[2]{%
                        {\raisebox{1pt}{$\scriptstyle\circ$}}\!%
                        \genfrac{}{}{0.5pt}{}{#1}{#2}%
                        \!\lower2pt\hbox{\rotatebox[origin=c]{-90}{{$\looparrowright$}}}}
\newcommand{\mgaction}[2]{%
                        {\raisebox{1pt}{$\scriptstyle :$}}\!%
                        \genfrac{}{}{0.5pt}{}{#1}{#2}%
                        \!\lower2pt\hbox{\rotatebox[origin=c]{-90}{{$\looparrowright$}}}}

\def\dgaction{\displaystyle\gaction}
\def\dggaction{\displaystyle\ggaction}

\def\smalll{\scriptsize}

\begin{document}

\title{Spinors, lattices, and classification \\  of integral Apollonian disk packings}

\author{Jerzy Kocik\\
{\small Department of Mathematics, Southern Illinois University, Carbondale, IL 62901}\\
{\small jkocik@siu.edu}} 

\date{}

\sloppy

\maketitle



\begin{abstract} 
\noindent
A parametrization of integral Descartes configurations (and effectively Apollonian disk packings)  
by pairs of two-dimensional integral vectors is presented.
The vectors, called here tangency spinors
defined for pairs of tangent disks, are spinors associated to the Clifford algebra 
for 3-dimensional Minkowski space.
A version with Pauli spinors is given.
The construction 
provides a novel interpretation to the known Diophantine equation 
parametrizing integral Apollonian packings.   
\\[5pt]
\end{abstract}
\noindent

~

\noindent
{\bf Notation:}  Throughout this paper both a disk and its curvature will be addressed by the same symbol.

\section{Introduction}

Integral Apollonian disk packings are remarkable geometric objects studied by number-theorists and geometers.
First of all, it might be surprising that such objects exist 
-- arrangements of tangent disks, all curvatures of which are integers.  
There are infinitely many of them and identifying them is a natural task.
The present paper shows that {\it any} pair of integral 2-vectors 
(effectively, a set of four integers $p$, $q$, $r$, $s$) 
generates an integral Descartes configuration and consequently an integral Apollonian packing:
$$
\mathbf a = \begin{bmatrix}  p\\ q  \end{bmatrix}, \quad  
\mathbf b = \begin{bmatrix}  r\\ s  \end{bmatrix}  
\qquad\Rightarrow\qquad
\begin{cases}
B_0 \  =  -|\mathbf a \times \mathbf b|\\
B_1 \  =  \ \ |\mathbf a \times \mathbf b|+\|\mathbf a\|^2 \\
B_2 \  =  \ \ |\mathbf a \times \mathbf b|+\|\mathbf b\|^2\\
B_3 \  =  \ \ |\mathbf a \times \mathbf b|+\|\mathbf a\pm\mathbf b\|^2
\end{cases}
$$ 
where $B_i$ denote the curvatures of the disks.
Formally, vectors $\mathbf a$ and $\mathbf b$ are spinors of 3-dimensional Minkowski space.
Here, they have a geometric context and are defined for any pair of tangent disks,
hence are called `tangency spinors''.  They are discussed in more detail in \cite{jk-spinor,jk-tessellation,jk-corona},
and reviewed in Section 2.
The above construction provides 
a novel interpretation to the known Diophantine equation that generates and parametrizes 
the integral packings \cite{jk-Diophantine}.
As a result, a one-to-one correspondence between irreducible integral Apollonian disk packings 
and certain  
irreducible integral sublattices of $\mathbb Z^2$ 
(each defined by a pair of spinors as the principal basis) emerges.
\def\smalll{\scriptsize}
%
$$
\begin{tikzpicture}[baseline=-0.8ex]
    \matrix (m) [ matrix of math nodes,
                         row sep=2.8em,
                         column sep=4em,
                         text height=3.8ex, text depth=3ex] 
   {
   \quad{\hbox{\sf A pair of} \atop \hbox{spinors}}  \quad  
      & \quad {\hbox{\sf Descartes} \atop \hbox{\sf configuration}} \quad  
      & \quad {\hbox{\sf Apollonian}\atop\hbox{\sf  packing }} \quad \\
      \quad {\hbox{\sf lattice}} \quad  & \quad ~ \qquad & \qquad \hbox{\sf corona} \quad   \\
     };
    \path[->]
        (m-1-1)   edge node[above] {$\hbox{ \sf\smalll  create}$} (m-1-2)
        (m-1-2)   edge  node[above] {$\hbox{ \sf\smalll  complete}$}(m-1-3)
        (m-2-3)   edge  node[above] {$\subset$} (m-2-1)
        (m-1-1) edge node[right]  {$\hbox{\sf\smalll generate}$} (m-2-1)
        (m-2-3) edge node[right] {$\subset$} (m-1-3)
        ;
\end{tikzpicture}   
$$

%
%

And now some basic notions organized in a form that should be easy to consult.
For more on the subject see \cite{Kate,GLM1,GLM4,jk-matrix,LMW,Mel,N,IS,Sod}.

\begin{itemize}[leftmargin=*]
\setlength\itemsep{0.2em}%

\vspace{-.05in}
\item
A {\bf disk} is an interior or the exterior of a circle.  
The two types of disks will be called {\bf inner} and {\bf outer}, respectively.
The former is assumed to have a positive curvature, the latter negative.

\vspace{-.05in}
\item
A {\bf tricycle} is a configuration of three mutually tangent disks, no pair overlapping.

\item
A {\bf Descartes configuration} is an arrangement of four mutually tangent circles.
Every tricycle may be completed to a Descartes configuration in two ways. 

\vspace{-.05in}
\item
{\bf Descartes formula} 
(found by Ren\'e Descartes in 1643) is a relation for Descartes configuration: 
\begin{equation}
\label{eq:Descartes}
2\, (A^2 + B^2 + C^2 + D^2) \ = \  (A + B + C + D)^2
\end{equation}
where $A$, $B$, $C$, and $D$ are the curvatures of the disks, i.e., reciprocals of their radii.  
Equation \eqref{eq:Descartes} solves {\bf Descartes problem}: 
given a tricycle, find curvature of fourth circle tangent to the three.
%
Due to quadratic nature of \eqref{eq:Descartes},
there are two solutions to the Descartes problem for the fourth disk:
$$
D=A+B+C\ \pm \ 2\sqrt{AB + BC + CD}
$$
We will call them {\bf conjugated} through disks $A,B,C$.
The resulting two Descartes configurations that differ by the fourth disk will also be called {\bf conjugated}. 



\vspace{-.05in}
\item
Descartes formula \eqref{eq:Descartes} may be viewed as a Diophantine equation, 
the integral solutions to which will be called {\bf Descartes quadruples}.
A Descartes  disk configuration is called {\bf integral} if the curvatures form a Descartes quadruple.
It is {\bf irreducible} if the curvatures are coprime, 
that is do not have a common factor, $\gcd(A,B,C,D)=1$. 
Here is an example of a pair of conjugated Descartes quadruples: 
$$
                   (\,2,\, 2, \, 3,\,  15\,)  \quad \hbox{and}\quad  (\,2,\,  2,\, 3,\, 35\,)
$$

\vspace{-.07in}
\item
A tricycle or Descartes configuration is {\bf everted}
if it contains a disk of negative curvature.

\vspace{-.05in}
\item
Every tricycle may be 
completed uniquely to an  
{\bf Apollonian disk packing} (called simply {\bf Apollonian packing}) 
by recursive inscribing new disks in the curvilinear triangular regions formed by the disks.
If an Apollonian packing contains an integral Descartes configuration then it is integral
i.e., its all disks have integral curvatures.
Consequently, the integral Apollonian disk packings may be classified by the integral Descartes 
configurations.

\vspace{-.05in}
\item
A Descartes configuration 
or tricycle is {\bf maximal} 
if  its completion to the Apollonian packing 
$A({\mathcal D})$ does not contain circles of smaller curvature than those in $\mathcal D$. 
This definition is motivated by the fact that an Apollonian packing 
may be grown from {\it any} Descartes quadruple it contains.
Reducing this redundancy to the maximal cases 
allows one to classify the integral Apollonian disk packings 
by maximal irreducible disk configurations.



\vspace{-.05in}
\item
The {\bf major disk}  in a disk configuration is the disk of negative or zero curvature if such exists.
Its boundary is called the {\bf major circle}. 
All integral Apollonian packings have a major circle 
(equivalent to the greatest among the circles).

\vspace{-.05in}
\item
The {\bf corona} is the set of all disks in the Apollonian packing that are tangent to a given disk. 
The {\bf major corona} is the corona at the major circle.
 For more on coronas see \cite{jk-corona}.

\end{itemize}

~

\begin{figure}[H] 
  \centering
  \includegraphics[width=2in,keepaspectratio]{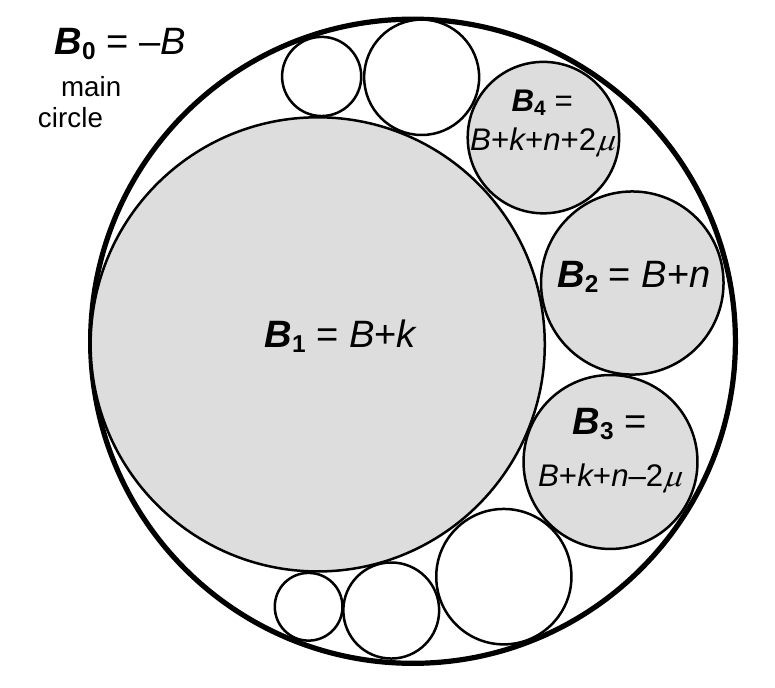}
  \caption{An Apollonian gasket and its greatest five  circles (smallest curvatures)}
  \label{fig:BBB}
\end{figure}

Let us also recall 
the theorem that defines an algorithm for generating all integral disk packings, 
copied here from \cite {jk-spinor}:

\begin{theorem} 
\label{thm:1}  \rm
[{\bf Parametrizing formula}]
{\sf
There is a one-to-one correspondence between the irreducible integral Apollonian gaskets 
and the irreducible quadruples of non-negative integers $B, \,k,\, n,\, \mu \in \mathbb{N}$ that satisfy 
\begin{equation}    
\label{eq:thm1}
  B^2+\mu^2 = kn
\end{equation}
with constraints   
\begin{equation}
\label{eq:conditions}
\begin{alignedat}{2}
  \hbox{(i)} &&\quad &0 \le \mu \le B/\surd 3, \\
  \hbox{(ii)} &&\quad &2\mu \le k \le n.
\end{alignedat}
\end{equation}
Every solution to \eqref{eq:thm1} corresponds to an integral Apollonian disk packing 
with the following quintet of the main curvatures:
\[
  (B_0,\,B_1\,,B_2,\,B_3,\,B_4)
	= (-B,\, B+k,\, B+n,\, B+k+n-2\mu,\, B+k+n+2\mu)
\]
The quadruples $(B_0,\,B_1\,,B_2,\,B_3)$ and $ (B_0,\,B_1\,,B_2,\,B_4)$
 are conjugated  \;  (see Figure \ref{fig:BBB}).
 }
 \end{theorem}

\begin{figure}[H]
  \centering
  \includegraphics[width=5.5in,keepaspectratio]{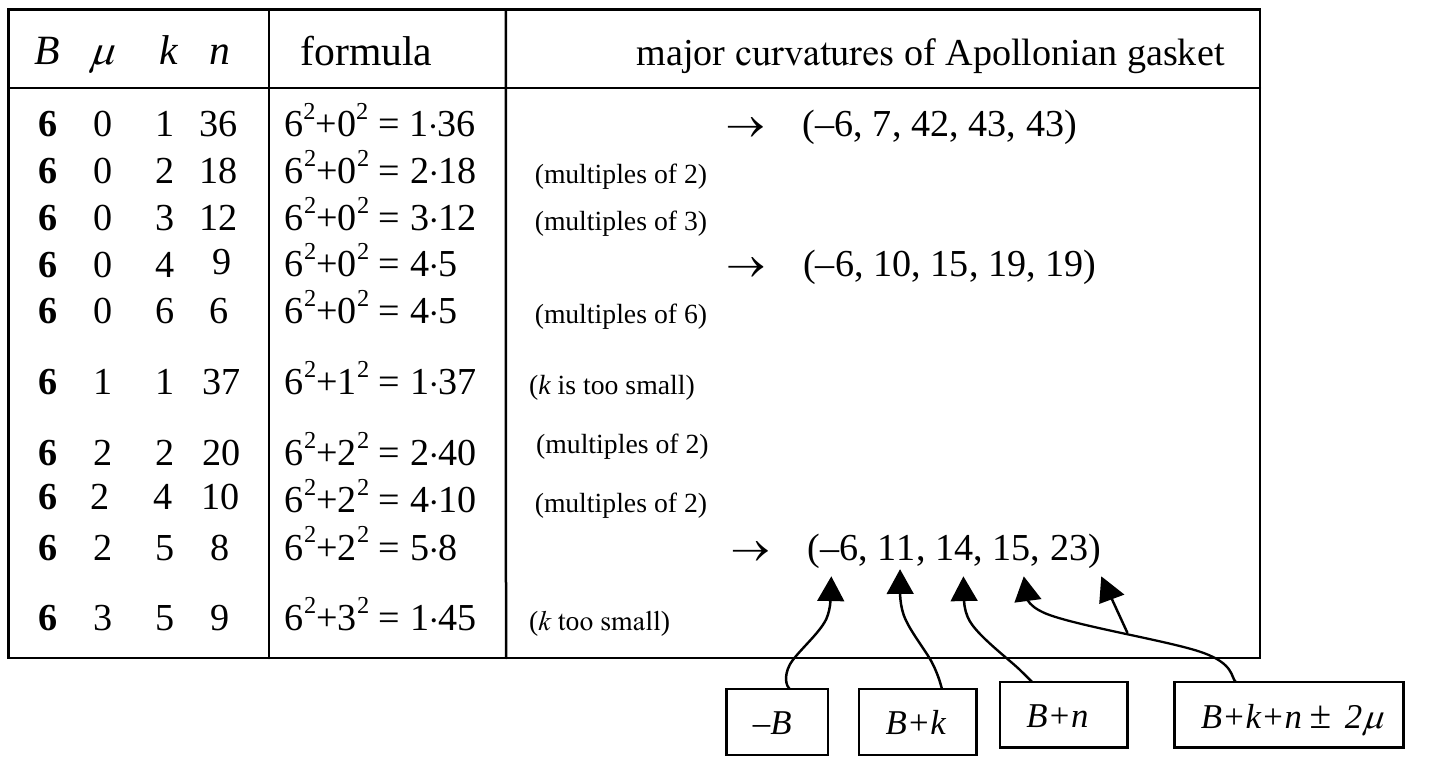}
  \caption{Producing the three maximal irreducible Descartes configurations for $B=6$} 
  \label{fig:produce}
\end{figure}

The derivation of this result was based on inversions of Apollonian strip 
with changing location of the point of tangency with the circle it was transformed.
Conditions \eqref{eq:thm1}--\eqref{eq:conditions} can easily be codified into an algorithm 
which lists all maximal irreducible integral Descartes quadruples in an organized fashion.  
Table in Figure \ref{fig:produce} shows the algorithm at work for $B=6$.

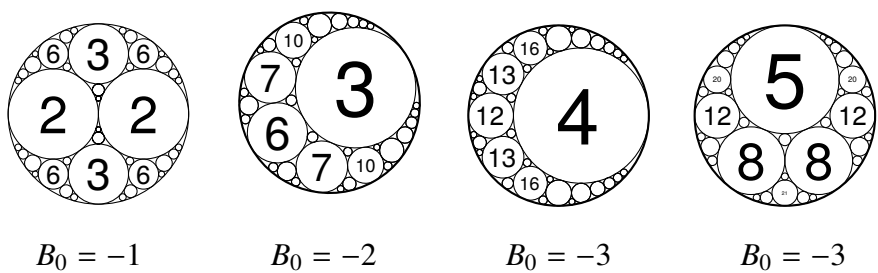
\begin{figure}[H]
\begin{center}

\begin{tikzpicture}[scale=1.2]  
\draw (0,0) circle (1);
\foreach \a/\b/\c/\s   in {
1 / 0 / 2 / 2,
-1 / 0 / 2 /2,
0 / 2 / 3   /1.7,
0 / -2 / 3   /1.7,
3 / 4 / 6 / 1,
-3 / 4 / 6 / 1,
3 / -4 / 6 / 1,
-3 / -4 / 6 / 1
}
\draw (\a/\c,\b/\c) circle (1/\c)
node  [scale=\s]  at (\a/\c,\b/\c)   {$\hbox{\sf \c}$}
;
\foreach \a /  \b / \c   in {
0/4/15, 0/6/35,  8/6/11, 5/12/14, 15/8/18, 
8/12/23, 7/24/26, 24/10/27, 21/20/30
}
\draw (\a/\c,\b/\c) circle (1/\c)          (-\a/\c,\b/\c) circle (1/\c)
          (\a/\c,-\b/\c) circle (1/\c)         (-\a/\c,-\b/\c) circle (1/\c)
;
\end{tikzpicture}
\hspace{.1in}
\begin{tikzpicture}[scale=2.4,rotate=30]
\draw [thick] (1/2,0) circle (1/2);
\foreach \a/\b/\c/\d/\s in {
2/ 0/ 3/ 1 /3, 
1/ 0/ 6/ 0/1.8, 
2/ 2/ 7/ 1/1.5, 
2/ -2/ 7/ 1/1.5, 
5/ 4/ 10/ 4/.7, 
5/ -4/ 10/ 4/.7
}
\draw (\a/\c,\b/\c) circle (1/\c)
        node [scale= \s]  at (\a/\c,\b/\c)  {$\hbox{\sf \c} \ $}
;
\foreach \a/\b/\c/\d in {
10/ 6/ 15/ 9,
10/ -6/ 15/ 9,
2/ 4/ 19/ 1,  2/ -4/ 19/ 1, 17/ 8/ 22/ 16, 17/ -8/ 22/ 16, 10/ 12/ 27/ 9, 10/ -12/ 27/ 9,
26/ 10/ 31/ 25, 26/ -10/ 31/ 25,  11/ 4/ 34/ 4, 11/ -4/ 34/ 4,  2/ 6/ 39/ 1,  2/ -6/ 39/ 1,
 5/ 12/ 42/ 4, 19/ 12/ 42/ 12, 37/ 12/ 42/ 36, 5/ -12/ 42/ 4, 19/ -12/ 42/ 12, 37/ -12/ 42/ 36,
26/ 20/ 43/ 25, 26/ -20/ 43/ 25, 17/ 24/ 54/ 16, 17/ -24/ 54/ 16, 50/ 14/ 55/ 49, 50/ -14/ 55/ 49,
35/ 20/ 58/ 28, 35/ -20/ 58/ 28
}
\draw (\a/\c,\b/\c) circle (1/\c)  
           (\a/\c,-\b/\c) circle (1/\c);
\end{tikzpicture}
\hspace{.1in}
%
%
%
\begin{tikzpicture}[scale=3.6]
\draw [thick] (2/3,0) circle (1/3);
\foreach \a/\b/\c/\d/\s in {
3/ 0/ 4/ 2  /3,
5/ 0/ 12/ 2/ 1,
6/ 2/ 13/ 3 /1,
6/ -2/ 13/ 3 /1,
9/ 4/ 16/ 6 /.7,
9/ -4/ 16/ 6 /.7 
}
\draw (\a/\c,\b/\c) circle (1/\c)
        node [scale= \s]  at (\a/\c,\b/\c)  {$\!\!\hbox{\sf \c}$}
;
\foreach \a/\b/\c/\d in {
14/ 6/ 21/ 11/2, 14/ -6/ 21/ 11/2, 21/ 8/ 28/ 18, 21/ -8/ 28/ 18, 30/ 10/ 37/ 27, 30/ -10/ 37/ 27,
15/ 4/ 40/ 6, 15/ -4/ 40/ 6, 23/ 12/ 48/ 14, 41/ 12/ 48/ 38, 23/ -12/ 48/ 14, 41/ -12/ 48/ 38,
30/ 4/ 61/ 15, 54/ 14/ 61/ 51, 30/ -4/ 61/ 15, 54/ -14/ 61/ 51, 39/ 20/ 64/ 30, 39/ -20/ 64/ 30,
38/ 12/ 69/ 23, 38/ -12/ 69/ 23, 69/ 16/ 76/ 66, 69/ -16/ 76/ 66, 
30/ 6/ 85/ 11, 54/ 20/ 85/ 39, 30/ -6/ 85/ 11, 54/ -20/ 85/ 39,
33/ 12/ 88/ 14, 63/ 28/ 88/ 54, 33/ -12/ 88/ 14, 63/ -28/ 88/ 54 
}
\draw (\a/\c,\b/\c) circle (1/\c)  
           (\a/\c,-\b/\c) circle (1/\c);
\end{tikzpicture}
\hspace{.15in}
%
%
%
\begin{tikzpicture}[scale=3.6]
\draw [thick] (0,2/3) circle (1/3);
\foreach \a/\b/\c/\d/\s in {
0/4/5/3/3, 
1/4/8/2/2, -1/4/8/2/2, 
-3/8/12/6/1,  3/8/12/6/1,  
-5/16/20/.3, 5/16/20/.3, 0/8/21/.2
}
\draw (\a/\c,\b/\c) circle (1/\c)
        node [scale= \s]  at (\a/\c,\b/\c)  {$\hbox{\sf \c}$}   ;
\foreach \a/\b/\c/\d in {
-8/ 16/ 29/ 11, 8/ 16/ 29/ 11, -7/ 28/ 32/ 26, 7/ 28/ 32/ 26,
 3/ 16/ 44/ 6, -3/ 16/ 44/ 6, 0/ 26/ 45/ 15, -9/ 44/ 48/ 42, 9/ 44/ 48/ 42,
-8/ 34/ 53/ 23,  -16/ 40/ 53/ 35, 8/ 34/ 53/ 23, 16/ 40/ 53/ 35, 
-15/ 28/ 56/ 18, 15/ 28/ 56/ 18, -21/ 40/ 68/ 30, -11/ 64/ 68/ 62,
21/ 40/ 68/ 30, 11/ 64/ 68/ 62, 0/ 34/ 77/ 15, 8/ 28/ 77/ 11, -8/ 28/ 77/ 11,
-16/ 58/ 77/ 47, 16/ 58/ 77/ 47, -13/ 88/ 92/ 86, 13/ 88/ 92/ 86
}
\draw (\a/\c,\b/\c) circle (1/\c)  
;
\end{tikzpicture}
\end{center}

\vspace{-.25in}

{\large
$$
B_0=-1
\hspace{.67in}
B_0=-2
\hspace{.67in}
B_0=-3
\hspace{.67in}
B_0=-3
$$
}

\vspace{-.2in}

\caption{First few Apollonian disk packings from Table 1}
\label{fig:examples}
\end{figure}

For any $B\in\mathbb N$ there are such irreducible maximal Descartes configuration 
and implied integral packings of the disk
of curvature $B_0=-B$   
(on average, about $B/3$ of them.)
An initial fragment of the resulting list is shown in Table \ref{fig:produce}.   
See also Figure \ref{fig:examples}. 
Note that, interestingly, columns $k$ and $n$ contain sums of squares.

~\\
{\bf Question:}  Is there a geometric picture behind the algebraic equation ?
\\\\
%
%
%
This is were the aforementioned ``tangency spinors''  come into play (Section 2).
Pairs of such spinors determine Descartes configurations (Section 3). 
With the help of the lattices that they define, a one-to-one parametrization of 
the integral packings is achieved (Section 4).
Consequently, the family of the irreducible integral Apollonian packings 
can be visualized as a dust of points in the celestial sphere, 
which matches with a similar image of the corresponding lattices in the hyperbolic upper half-plane
(Section 5).

{\small
$$\left\{
\begin{array}{ll}
\hbox{\sf irreducible}\\
\hbox{\sf integral}\\
\hbox{\sf Apollonian}\\
\hbox{\sf disk}\\
\hbox{\sf packings}
\end{array}
\right.
\ 
\xleftrightarrow{\hbox{1-to-1}}
\ 
\left\{
\begin{array}{ll}
\hbox{\sf maximal}\\
\hbox{\sf irreduciblel}\\
\hbox{\sf integral}\\
\hbox{\sf Descartes}\\
\hbox{\sf configurations}
\end{array}
\right.
\ 
\xleftrightarrow{\hbox{1-to-1}}
\ 
\left\{
\begin{array}{ll}
\hbox{\sf reduced}\\
\hbox{\sf integral}\\
\hbox{\sf spinors}\\
\hbox{\sf (2D vectors)}
\end{array}
\right.
\xleftrightarrow{\hbox{1-to-1}}
\ 
\left\{
\begin{array}{ll}
\hbox{\sf classes of}\\
\hbox{\sf irreducible}\\
\hbox{\sf sublattices}\\
\hbox{\sf of } \mathbb Z^2
\end{array}
\right.
\ 
\xleftrightarrow{\hbox{1-to-1}}
\ 
\left\{
\begin{array}{ll}
\hbox{\sf dust in the}\\
\hbox{\sf fundamental}\\
\hbox{\sf domain}\\
\hbox{\sf of Dedekind}\\
\hbox{\sf tessellation}
\end{array}
\right.
$$
}

\newpage

\begin{table}[htbp] %
  \centering
  \includegraphics[scale=.8]{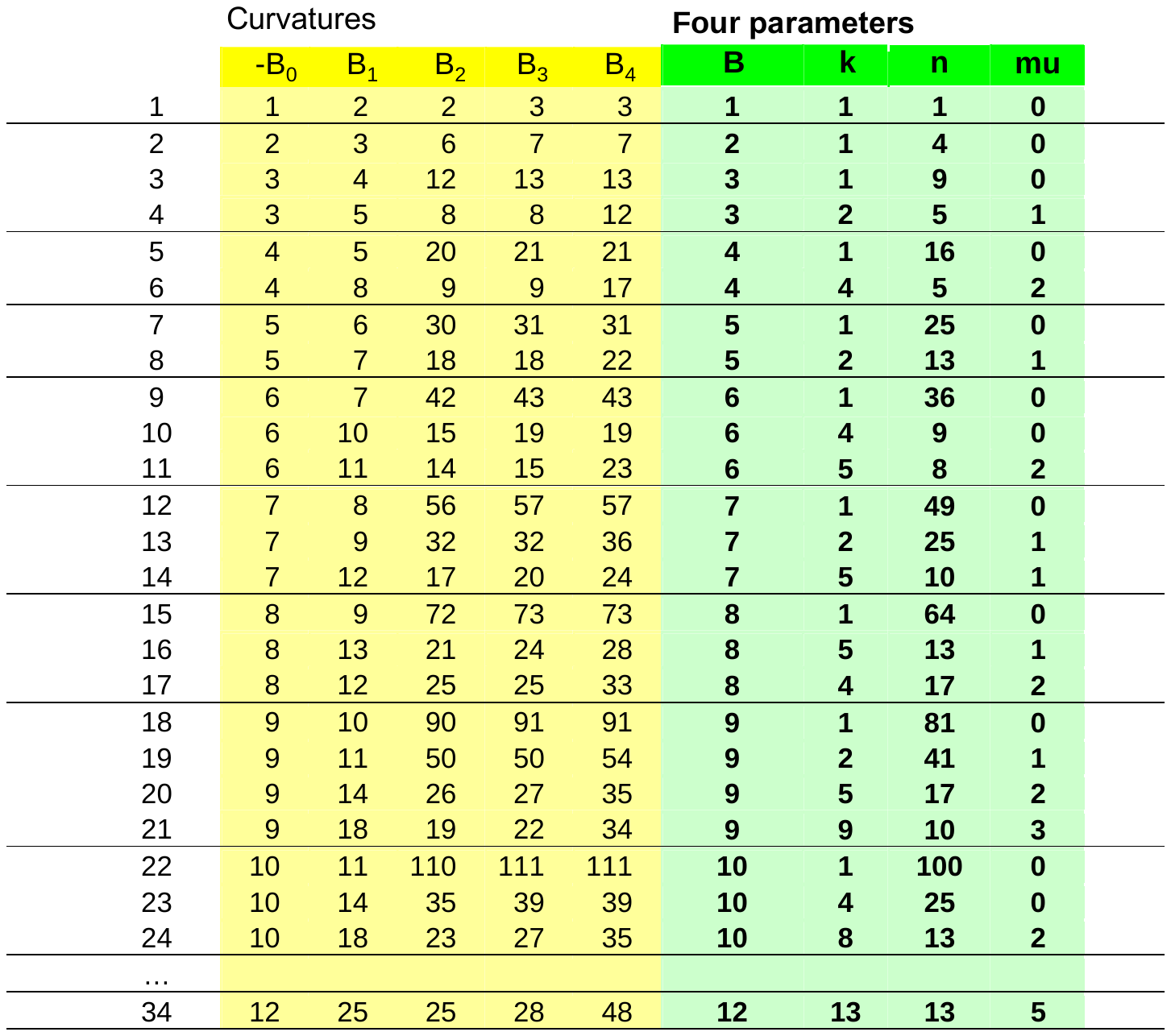}  
  \caption{Bend quintets for Apollonian integer disk packings for principal curvatures 1 through 10, 
  including four parameters.}
  \label{fig:results}
\end{table}


We finish the paper with a surprising coda.
The pair of tangency spinors define
a Pauli spinor, a vector in $\mathbb C^2$, 
well-known in theoretical physics.
The parametrizing formula obtains a form of 
vanishing determinant of the Hermitian matrix 
defined by the Kronecker power of this spinor:
$$
\begin{bmatrix}
   k & \mu -iB\\
   \mu +iB & n \end{bmatrix}
   \ = \ 
\begin{bmatrix} \mathbf a  \\ \mathbf b \end{bmatrix}
\otimes
\begin{bmatrix} \mathbf a \\ \mathbf b \end{bmatrix}^*
\qquad  \xrightarrow{\quad\rm det \quad} \ 0 \,,  
$$
where the star denotes Hermitian conjugation. 
This brings the world of theoretical physics with the (1+3) space-time and 
associated spinors into the picture.

\newpage
\section{Spinor space for tangent disks}

In his section we review briefly the notion of {\bf tangency spinor} of disks.
For details,  proofs, and motivation  see \cite{jk-spinor}.

\subsection{Spinor space}

In the following, by the spinor space we will mean 
a two dimensional real vector space, 
which may be also emulated with 1-dimensional complex space (Argand plane), 
$\mathbb R^2 \cong \mathbb C$.
Typical vectors are:
$$
\mathbf a = \begin{bmatrix}  x\\ y  \end{bmatrix} 
\qquad
\mathbf b = \begin{bmatrix}  x' \\ y'  \end{bmatrix}  
\qquad
$$
The space is equipped with two structures,
the Euclidean inner product (``dot product')' and the symplectic product (``cross-product''), 
both with values in real numbers:
$$
\begin{array}{rclll} 
\hbox{inner product:}        &\mathbf a,\mathbf b  \quad \mapsto\quad  &\mathbf a\bigcdot \mathbf b =  xx' + yy'    \\[5pt]
\hbox{symplectic product:} &\mathbf a,\mathbf b  \quad \mapsto\quad  &\mathbf a\times \mathbf b =  xy' - x'y  
           \ \equiv\det [\mathbf a\,\mathbf b]    \\[5pt]
\end{array}
$$
We also define `'symplectic conjugation''
$$
\mathbf a^{\symp} \ = \  
         \begin{bmatrix}  x\\ y  \end{bmatrix}^{\symp} \ = \  
         \begin{bmatrix}  -y\\ \phantom{-}x  \end{bmatrix} \ = \ 
         \begin{bmatrix}  0&\!\!\!\!-\!1\\ 1&0  \end{bmatrix}\begin{bmatrix}  x\\ y  \end{bmatrix}
$$
The two structures are related:
$$
\mathbf a\times \mathbf b \ = \  \mathbf a^{\symp}\bigcdot \mathbf b  
$$
Other identities are readily implied:
$$
 (\mathbf a^{\symp})^{\symp} = - \mathbf a\,, 
\quad \mathbf a\bigcdot \mathbf b = \mathbf a\times \mathbf b^{\symp}\,, 
\quad \mathbf a^{\symp} \bigcdot \mathbf b^{\symp} = \mathbf a\bigcdot \mathbf b \,,
\quad \mathbf a^{\symp} \times \mathbf b^{\symp} = \mathbf a\times \mathbf b \,,
\quad \mathbf a^{\symp} \times \mathbf b = \mathbf b^{\symp} \times \mathbf a \,.
$$
The squares are:
$$
\mathbf a\bigcdot \mathbf a =\|\mathbf a\|^2  =x^2 + y^2\,,
\qquad  \mathbf a\times \mathbf a = 0\,,
\qquad \mathbf a\bigcdot \mathbf a^{\symp} = 0\,.
$$
{\bf Interpretation via complex numbers.}
When the spinor space is represented by complex numbers,
$$
a = \begin{bmatrix}  x\\ y  \end{bmatrix}  \ = \ x+yi
\qquad\hbox{and}\qquad
b = \begin{bmatrix}  x' \\ y'  \end{bmatrix}  \ = \ x'+y' i
$$
then the above structures are expressed as follows:
$$
\begin{array}{rrlll} 
\hbox{inner product:}           &  a,b  \quad \mapsto\quad   &a\bigcdot b = \frac{1}{2} (\bar a b   + a \bar b)     \\[5pt]
\hbox{symplectic product:}   &  a,b  \quad \mapsto\quad   &a\times b = \frac{1}{2i} (\bar a b   - a \bar b)     \\[5pt]
\hbox{conjugation:}             &  a \quad \mapsto \quad  &a^{\symp} = ai \\[5pt]
\end{array}
$$
Note that by ``conjugation'' we mean ``symplectic conjugation'' (denoted by star $\symp$).
Not to be confused with ``complex conjugation'',
always called by its full name and denoted by the regular asterisk $\ast$.

\subsection{Spinors and Descartes}

\noindent
{\bf Definition:}  Let $A$ and $B$ be an ordered pair of mutually tangent disks 
of radii  $r_A$ and $r_B$ and centered at $C_A$ and $C_B$, respectively, 
in a plane identified with complex numbers, $\mathbb C\cong \mathbb R^2$.
Interpret the vector joining the centers as a complex number $z= (C_AC_B)$.
The {\bf tangency spinor} of the two disks is  a complex number $u$
or equivalently as a 2-vector $\mathbf u$
defined as
$$
u \  = \ \pm \sqrt{\frac{z}{r_Ar_B}} 
\in \mathbb C 
\qquad
\mathbf u \  = \ \pm  
\begin{bmatrix}\hbox{\rm Re}\, u\\ \hbox{\rm Im} \, u \end{bmatrix}   
\in \mathbb R^2
$$
The spinor is defined up to a sign since $(-u)^2 = u^2$.
Also, the spinor depends on the order of disks:  
 if $u$ is a spinor for $(AB)$,  then the spinor for  $(BA)$  is $u^{\symp} = iu$  
 (again, both up to sign).
\\

In graphical representation we will mark a spinor by an arrow that indicates 
the order of circles, and will label it by its matrix value.
The geometric interpretation and motivation can be found in \cite{jk-spinor} and in Appendix A. 
The proofs of the main properties listed below are in \cite{jk-spinor}.
For economy, capital letters will denote both circles and their curvatures.

\begin{figure}[h]
\centering
\hrule

~\\[7pt]

\begin{tikzpicture}[scale=1.2]
\draw[fill=gold!10]  (-.19, 0) circle (.7);
\draw[fill=gold!10]  (1, 0.5) circle (.58);
\draw [->, ultra thick, color=black] (.25, .12) -- (.75, .37);
\node  [scale=1.1]  at (-.5, .2) {\sf  A}; 
\node  [scale=1.1]  at (1.2, .7) {\sf  B}; 
\node  [scale=1.1]  at (.1, -.1) {\sf  u}; 
\end{tikzpicture}
\qquad\qquad
\begin{tikzpicture}[scale=1.1]
\clip (-.95, .3) rectangle  (1.6,2.45);
\draw[fill=gold!10]  (0,0) circle (1);
\draw[fill=gold!10]  (-.2, 1.7) circle (.7);
\draw[fill=gold!10]  (1, 1.226) circle (.58);
\draw [->, ultra thick, color=black] (-.115, .75) -- (-.15,1.25);
\draw [->, ultra thick, color=black] (.53, .62) -- (.83,1.0);
\node  [scale=1.1]  at (-.6, .43) {\sf  C}; 
\node  [scale=1.1]  at (-.5, 1.9) {\sf  A}; 
\node  [scale=1.2]  at (1.2, 1.45) {\sf  B}; 
\node  [scale=1.1]  at (-.1, .5) {\sf  a}; 
\node  [scale=1.1]  at (.45, .48) {\sf  b}; 
\end{tikzpicture}
\qquad\qquad
\begin{tikzpicture}[scale=1.1]
\clip (-.95, .3) rectangle  (1.6,2.45);
\draw[fill=gold!10]  (0,0) circle (1);
\draw[fill=gold!10]  (-.2, 1.7) circle (.7);
\draw[fill=gold!10]  (1, 1.226) circle (.58);
\draw[thick, dotted]  (.3, 1.05) circle (.41);
\draw [->, ultra thick, color=black] (-.115, .75) -- (-.15,1.25);
\draw [->, ultra thick, color=black] (.53, .62) -- (.83,1.0);
\node  [scale=1.2]  at (-.1, .5) {\sf  a}; 
\node  [scale=1.2]  at (.45, .48) {\sf  b}; 
\node  [scale=1]  at (.2, 1.6) {\sf T}; 
\end{tikzpicture}


\begin{tikzpicture}[scale=1.1]
\node  [scale=1]  at (-4.5,0) {$|\mathbf u|^2 = A+B $}; 
\node  [scale=1]  at (0, 0) {$\mathbf a\times \mathbf b =\pm \, C$}; 
\node  [scale=1]  at (3.7, 0) {$\mathbf a\bigcdot \mathbf b = T$};
\node  [scale=1]  at (4.4, 0) {~};
\node  [scale=1]  at (0,.3) {~};
\end{tikzpicture}


~

\begin{tikzpicture}[scale=1.1]
\clip (-.95, .3) rectangle  (1.6,2.45);

\draw[fill=gold!10]  (0,0) circle (1);
\draw[fill=gold!10]  (-.2, 1.7) circle (.7);
\draw[fill=gold!10]  (1, 1.226) circle (.58);
\draw [->, ultra thick, color=black] (-.115, .75) -- (-.15,1.25);
\draw [<-, ultra thick, color=black] (.53, .62) -- (.83,1.0);
\draw [->, ultra thick, color=black] (.25, 1.54) -- (.75,1.33);
\node  [scale=1.1]  at (-.1, .5) {\sf  a}; 
\node  [scale=1.1]  at (1, .98) {\sf  c}; 
\node  [scale=1.1]  at (.1, 1.7) {\sf  b}; 
\end{tikzpicture}
\qquad\qquad
\begin{tikzpicture}[scale=1.1]
\clip (-.95, .3) rectangle  (1.6,2.45);
\draw[fill=gold!10]  (0,0) circle (1);
\draw[fill=gold!10]  (-.2, 1.7) circle (.7);
\draw[fill=gold!10]  (1, 1.226) circle (.58);
\draw [fill=gold!10] (.32, 1.07) circle (.1167);
\draw [->, ultra thick, color=black] (.115, .6) -- (.24,.9);
\node  [scale=1.1]  at (-.1, .5) {\sf  a}; 
\draw [->, ultra thick, color=black] (.83, 1.22) -- (.53,1.12);
\node  [scale=1.1]  at (1, 1.2) {\sf  c}; 
\draw [->, ultra thick, color=black] (.03, 1.45) -- (.23,1.2);
\node  [scale=1.1]  at (0, 1.7) {\sf  b}; 
\end{tikzpicture}
\qquad\qquad
\begin{tikzpicture}[scale=1.1]
\clip (-3, -.08) rectangle  (0.1, 2.2);
\draw[fill=gold!10]  (0,-5) circle (6);
\draw[fill=gold!10]  (-2, 1.4) circle (.7);
\draw[fill=gold!10]  (-.8, 1.455) circle (.5);
\draw[fill=gold!10]  (-1.257, 1.01) circle (.136);
\draw [->, thick, color=black] (-1.6, .22) -- (-1.75, .67);
\node  [scale=1.1]  at (-1.55, .01) {\sf  a}; 
\draw [->, thick, color=black] (-.65, .37) -- (-.73, .89);
\node  [scale=1.1]  at (-.58, .14) {\sf  b}; 
\draw [->, thick, color=black] (-1.12, .51) -- (-1.2, .82);
\node  [scale=1]  at (-1.1, .28) {\sf  a+b}; 
\end{tikzpicture}


\begin{tikzpicture}[scale=1.2]
\node  [scale=1]  at (-4.5,0) {$\mathbf a+\mathbf b+\mathbf c \ =\ 0$}; 
\node  [scale=1]  at (-.5, 0) {$\mathbf a+\mathbf b+\mathbf c \ =\ 0$}; 
\node  [scale=1]  at (3.5, 0) {$\mathbf a+\mathbf b \ =\ \mathbf c$};
\node  [scale=1]  at (4, 0) {~};
\node  [scale=1]  at (0,.3) {~};
\node  [scale=1]  at (-4.5,-.5) {``curl\,$u \ =\ 0$''}; 
\node  [scale=1]  at (-.5, .-.5) {``div\,$u\ =\ 0$''};
\end{tikzpicture}

~

\hrule

~

\caption{Spinors and disks.}
\label{fig:6figs}
\end{figure}
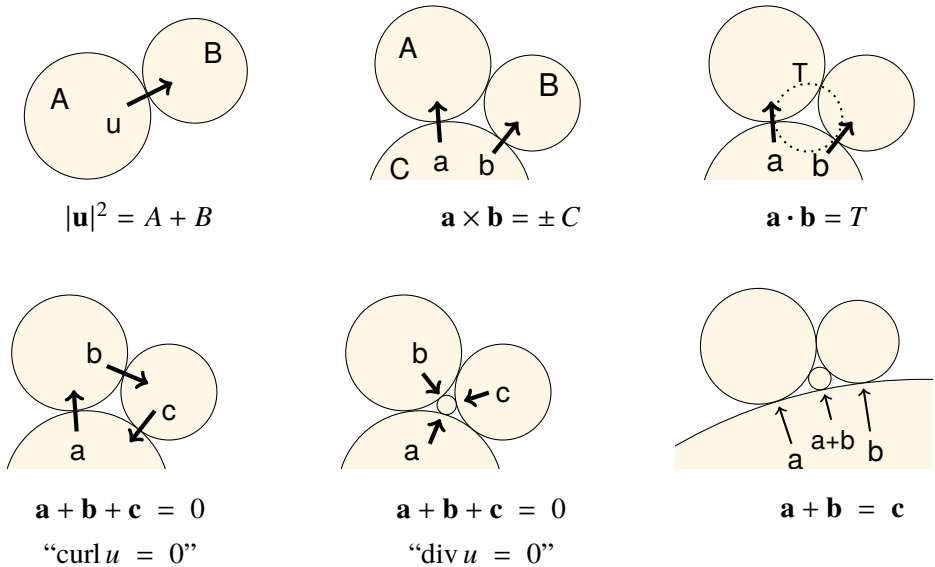


%
%

\noindent
{\bf Theorem 2.1.} 
\label{thm:curvs}
{\sf
 If $u$ is the tangency spinor for two tangent disks of curvatures $A$ and $B$, respectively,   
(Figure \ref{fig:6figs}, top left)  then 
\begin{equation}
\label{eq:thm3}
                     |\mathbf u|^2 = A + B
\end{equation}
%
}
\noindent
{\bf Theorem 2.2.}  [{\bf curvatures from spinors}]    
\label{thm:curv}
{\sf
In the system of three mutually tangent circles, the symplectic product of two spinors directed outward from 
(respectively inward into) one of the circles is equal (up to sign) to its curvature. 
Following notation of Figure \ref{fig:6figs} top center:
\begin{equation}
\label{eq:thm4}
          \pm \, C  \ = \ \mathbf a\times \mathbf b \   =   \ \det [\mathbf a|\mathbf b]
\end{equation}
}
\noindent
{\bf Theorem 2.3.}  [{\bf mid-circle from spinors}]    
\label{thm:mid}
{\sf 
For the situation as above, the dot product of the spinors equals the curvature
of the mid-circle $(ABC)$, i.e., the circle through the points of their tangency. 
Following the notation of Figure \ref{fig:6figs}, top right:
\begin{equation}
\label{eq:thm42}
           T  =  \mathbf a\bigcdot \mathbf b \qquad (\hbox{respectively, } T = - \mathbf a\bigcdot \mathbf b)
\end{equation}
}
\noindent
{\bf Theorem 2.4.}  [{\bf spinor curl}]. 
\label{thm:curl}
{\sf
The signs of the three tangency spinors between three mutually tangent circles
(Figure \ref{fig:6figs}, left)  
may be chosen so that 
\begin{equation}
\label{eq:curl}
\mathbf a  +  \mathbf b  +  \mathbf c  =  0  \qquad\quad         [ \hbox{\rm ``curl}\ \mathbf u = 0 \hbox{''}]   
\end{equation}
%
}
\noindent
{\bf Theorem 2.5.}
Let $A$, $B$, $C$, and $D$ be four circles in a Descartes configuration.  
\\[5pt]
{\sf 
{\bf [A. Vanishing divergence]:}  
If $a$, $b$ and $c$ are tangency spinors for pairs $AD$, $BD$ and $CD$  
(see Figure \ref{fig:6figs}, bottom center),
then their signs may be chosen so that 
\begin{equation}
\label{eq:thm5a}
 			\mathbf a + \mathbf b + \mathbf c = 0 	\qquad	[\hbox{``\rm div}\, \mathbf u = 0\hbox{''}]                                   
\end{equation}
The same property holds for the outward oriented spinors.
\\[5pt]
{\bf [B. Additivity]:}  If $a$ and $b$ are spinors of tangency for pairs $CA$ and $CB$ 
(see  Figure \ref{fig:6figs}, bottom right), 
then there is a choice of signs so that the sum 
\begin{equation}
\label{eq:thm5b}
                                     \mathbf  c = \mathbf a + \mathbf b
\end{equation}
is the tangency  spinor of $CD$.  
}

%

~

We shall need also the following concept:
Two spinors from disk $A$ and $B$ to disk $C$ are {\bf adjacent} if $A$, $B$, and $C$ are tangent
and there exists disk $D$ such that the whole quadruple forms a Descartes configuration.
The same goes for reciprocal spinors, i.e.,  {\it to} $C$.  
For instance, spinors $\mathbf a$ and $\mathbf b$ in Figure \ref{fig:nice} are adjacent, 
while $[3,0]^T$ and $[-1,-4]^T$ are not.

~\\
{\bf Remark 1:}
Theorem 2.5 may be viewed as a ``square root of Descartes Theorem'', 
since squaring the formula in a careful way reproduces Descartes formula \cite{jk-spinor}.
(The two version of this theorem are equivalent and differ by a sign of one spinor and geometric interpretation).
\\[7pt]
{\bf Remark 2:}
It should be noted that the ``divergence'' and the ``curl'' theorems have local character.
Extending the spinor vector field to the Apollonian disk packing cannot be done consistently
due to topological obstruction, quite like the non-existence of smooth non-vanishing vector field on a sphere. 
\\[7pt]
{\bf Remark 3:} the name tangency ``spinor'' is justified by the fact that
it represents spinor for the (1+2)-dimensional Minkowski space 
in which Pythagorean triples are null-vectors (metaphorically, photons)
\cite{jk-Clifford}.
Consequently, they behave like electron spins  known in quantum physics:
rotating a configuration of circles around by $360^{\circ}$ causes the spinors rotate by only $180^{\circ}$,
i.e., they change the signs.
It takes two full rotations of a circle configuration to bring the spinors to their original signs. 
\\[7pt]
{\bf Notation:} A spinor from disk $A$ to $B$ will be denoted sometimes as
$$
\spin(A,B)
$$
Clearly,  $\spin(B,A)=\pm (\spin(A,B))^{\symp}$ 
(or $\spin(B,A)=\pm i\, (\spin(A,B))$ when represented as complex numbers).


\section{Integral Descartes configurations from spinors}

\begin{figure}[H]
\centering
\begin{tikzpicture}[scale=12,rotate=60]
\draw [fill=yellow!30, rotate=-60](-.66,.53) rectangle (0,1.03) ;
\draw [fill=white,thick] (3/6,4/6) circle (1/6);
\draw [fill=yellow!30, thick] (6/11,8/11) circle (1/11);
\draw [fill=yellow!30, thick] (7/14,8/14) circle (1/14);
\draw [fill=yellow!10, thick] (6/15,10/15) circle (1/15);
\draw [fill=yellow!10, thick] (14/23, 14/ 23) circle (1/23);
\draw [thick] (6/11,8/11) circle (1/11)
                   node [scale= 1.8]  at (6/11,8/11)  {$B_1$};
\draw [thick] (7/14,8/14) circle (1/14)
                   node [scale= 1.8]  at (7/14,8/14)  {$B_2$};
\draw [thin, opacity=.35] (6/15,10/15) circle (1/15)
                   node [scale= 1.8]  at (6/15,10/15)  {$B_3$};
\draw [thin, opacity=.35] (14/23, 14/ 23) circle (1/23)
                   node [scale= 1.8]  at (14/23, 14/ 23)  {$B_4$};
\node [scale= 1.7]  at (.34, .83)  {$B_0$};
\draw [black,->, line width=1.2mm] (3/5+.012+.0096, 4/5+.015+.0120)--(3/5-.012-.012, 4/5-.015-.015);
\draw [black,->, line width=1.2mm] (4/8, 4/8-.040)--(4/8, 4/8+.035);
\draw [blue,->, line width=.57mm] (3/9-.03, 6/9)--(3/9+.02, 6/9);
\draw [blue,->, line width=.57mm] (11/17+.02, 10/17-.01)--(11/17-.02, 10/17+.01);

\draw node  [scale=1.4]  at (3/5+.042, 4/5+.045) {\color{black} $\mathbf a$};
\draw node  [scale=1.4]  at (4/8-.0, 4/8-.08) {\color{black} $\mathbf b$};
\draw node  [scale=1.2]  at (3/9-.049-.03, 6/9+.027) {\color{black} $\mathbf a \!-\! \mathbf b$};
\draw node  [scale=1.2]  at (11/17+.049+.01, 10/17-.02) {\color{black} $\mathbf a\! +\! \mathbf b$};

\end{tikzpicture}
\caption{Everted Descartes configuration and two spinors}
\label{fig:Descartes}
\end{figure}
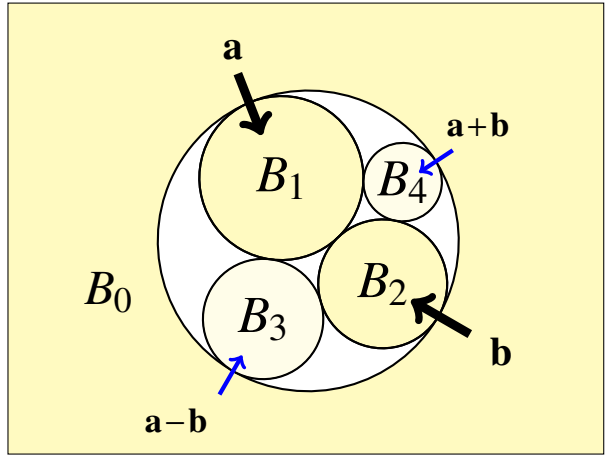

\noindent
{\bf Proposition 3.1:}
{\sf
Suppose two tangent disks of positive curvatures $B_1$ and $B_2$ are inscribed 
in a disk of negative curvature $B_0$ (see Figure \ref{fig:Descartes}). 
Suppose the tangency spinors are   
$$
\mathbf a = \spin(B_0,B_1)\,,
\qquad
\mathbf b = \spin(B_0,B_2)\,,
\qquad M=[\mathbf a \mathbf b].  
$$
Denote $B\phantom{_0} =  |\mathbf a \times \mathbf b |  = |\/\det M\/|$. 
Then the following are true
\begin{equation}
\label{eq:BBB1}
\begin{array}{rlc}
B_0 & =\ -B     \\[3pt]
B_1 & =\ B+ \|\mathbf a \|^2     \\[3pt]   
B_2 & =\ B+ \|\mathbf b \|^2     \\[3pt]   
\end{array}
\end{equation}
Moreover, 
the following quadruples 
satisfy Descartes formula 
$$
(B_0,B_1,B_2,B_3)  \qquad \hbox{and} \qquad (B_0,B_1,B_2,B_4)\,,
$$
where 
$$
B_{3,4}  =\ B +  \|\mathbf a\pm \mathbf b \|^2     =
   \ B +  \|\mathbf a\|^2+ \|\mathbf b \|^2  \pm \mathbf a\bigcdot\mathbf b  
$$
}

\noindent
{\bf Proof:}  
Direct consequence of Theorems 2.1--2.5 of the previous section.  Indeed, one has
$$
\begin{array}{rlll}
|\mathbf a \times \mathbf b | &= B\phantom{_0} = |\/\det M\/| &\quad& \hbox{\small (by Thm 2.2)}  \\[3pt]   
\|\mathbf a \|^2& =\ B_1 + B_0  && \hbox{\small (by Thm 2.1)}  \\[3pt]   
\|\mathbf b \|^2& =\ B_2 + B_0  && \hbox{\small (by Thm 2.1)}  \\[3pt]   
\|\mathbf a+\mathbf b \|^2& =\ B_3 + B_0   && \hbox{\small (by Thm 2.5B)}  \\[3pt]   
\end{array}
$$
Denote $B_0=-B$ and invert the equations to get the claim.
\qed

~

Note that the scalar $B$ may be defined by any pair of the above spinors:
$$
B =  |\mathbf a \times \mathbf b | =  |\mathbf a \times \mathbf c | =  |\mathbf b \times \mathbf c | 
$$
where $ \mathbf c  =  \mathbf a \pm \mathbf b$.
Indeed, $\mathbf a \times \mathbf c = \mathbf a \times (\mathbf a +\mathbf b) = \mathbf a \times \mathbf b$,
and similarly for the second equation.  

~

\begin{figure}
\centering
\begin{tikzpicture}[scale=20]

\draw [thick] (3/6,4/6) circle (1/6);
\draw [fill=yellow!30, thick] (6/11,8/11) circle (1/11);
\draw [fill=yellow!30, thick] (7/14,8/14) circle (1/14);
\draw [fill=yellow!10, thick] (6/15,10/15) circle (1/15);
\foreach \a/\b/\c/\d/\s in {
 6/ 8/ 11/ 94/1.7,
 7/ 8/ 14/ 8/1.7,
 6/ 10/ 15/ 9/1.7,
14/ 14/ 23/ 17/1,
11/ 20/ 26/ 20/1,
14/ 20/ 35/ 17/.9,
27/ 28/ 42/ 36/.7,
22/ 38/ 47/ 41/.5,
30/ 28/ 51/ 33/.5,
22/ 44/ 59/ 41/.5
}
\draw (\a/\c,\b/\c) circle (1/\c)
        node [scale= 1.7*\s]  at (\a/\c,\b/\c)  {\sf \c \; }
;

\foreach \a/\b/\c/\d in {
30/ 38/ 71/ 33,
46/ 50/ 71/ 65,
27/ 44/ 74/ 36,
39/ 64/ 78/ 72,
41/ 56/ 86/ 56,
54/ 50/ 95/ 57,
57/ 64/ 102/ 72,
70/ 68/ 107/ 89,
49/ 80/ 110/ 80,
39/ 80/ 110/ 72,
71/ 80/ 110/ 104,
62/ 98/ 119/ 113,
75/ 68/ 122/ 84,
54/ 64/ 123/ 57,
54/ 100/ 123/ 105,
46/ 80/ 131/ 65,
57/ 80/ 134/ 72,
54/ 110/ 143/ 105,
86/ 80/ 155/ 89,
97/ 104/ 158/ 128,
102/ 118/ 159/ 153,
91/ 140/ 170/ 164,
81/ 136/ 174/ 144,
62/ 128/ 179/ 113,
105/ 104/ 182/ 120,
75/ 100/ 186/ 84,
86/ 98/ 191/ 89,
70/ 110/ 191/ 89
}
\draw (\a/\c,\b/\c) circle (1/\c)  
node [scale= 45/\c]  at (\a/\c,\b/\c)  {$\c$}
;


\node [scale= 1.7]  at (.34, .83)  {$\frac{\hbox{\sf -3,-4}}{\hbox{-6}}$};

\draw [red,->, line width=1.2mm] (3/5+.012+.008, 4/5+.015+.01)--(3/5-.012-.004, 4/5-.015-.005);
\draw [red,->, line width=1.2mm] (4/8, 4/8-.032)--(4/8, 4/8+.028);
\draw [blue,->, line width=.7mm] (3/9-.03, 6/9)--(3/9+.02, 6/9);
\draw [blue,->, line width=.7mm] (11/17+.02, 10/17-.01)--(11/17-.02, 10/17+.01);

\draw node  [scale=1.2]  at (3/5+.042, 4/5+.045) {\color{black} $\mathbf a = \begin{bmatrix}\;\mathbf 1\\ \hbox{-}\mathbf 2\end{bmatrix}$};
\draw node  [scale=1.2]  at (4/8-.07, 4/8-.043) {\color{black} $\mathbf b = \begin{bmatrix}\mathbf 2 \\ \mathbf 2\end{bmatrix}$};
\draw node  [scale=1]  at (3/9-.049-.05, 6/9) {\color{black} $\mathbf a+\mathbf b=\begin{bmatrix}3 \\ 0\end{bmatrix}$};
\draw node  [scale=1]  at (11/17+.049+.05, 10/17-.02) {\color{black} $\begin{bmatrix}\hbox{-}1 \\ \hbox{-}4\end{bmatrix}=\mathbf a-\mathbf b$};

%
%
\end{tikzpicture}
\caption{Packing $(-6,11,14,15)$ and the principal spinors}
\label{fig:nice}
\end{figure}
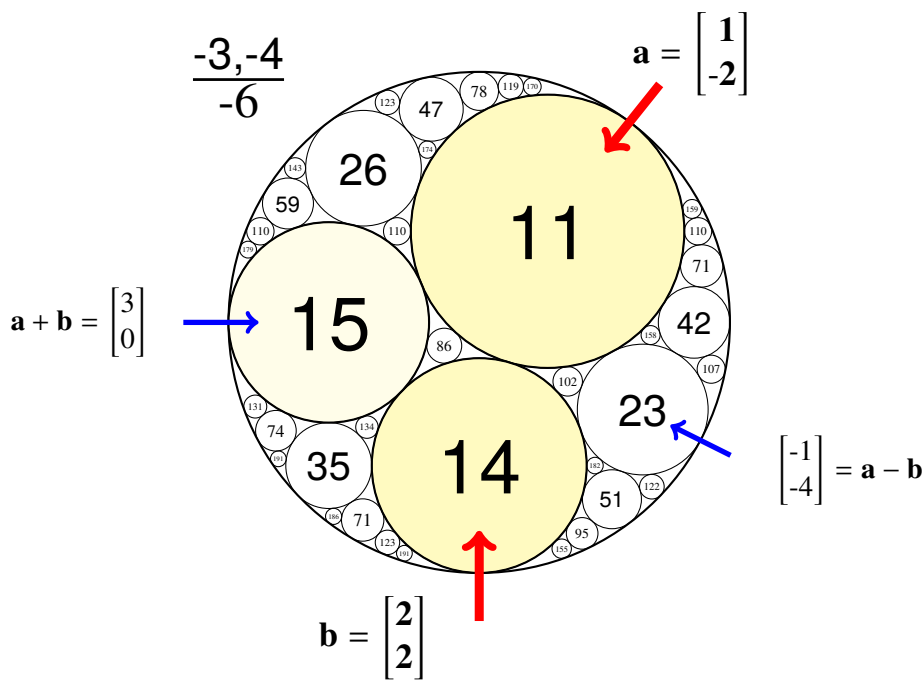

\noindent
{\bf Proposition 3.2:} 
{\sf Any two integral vectors $\mathbf a$ and $\mathbf b$ determine an everted Descartes configuration
with integral curvatures $(B_0, B_1, B_2, B_3)$ defined
\begin{equation}
\label{eq:main}
\begin{array}{llc}
B_0 & =\ -|\mathbf a \times \mathbf b |      \\[3pt]
B_1 & =\ \ |\mathbf a \times \mathbf b | + \|\mathbf a \|^2     \\[3pt]   
B_2 & =\ \ |\mathbf a \times \mathbf b | + \|\mathbf b \|^2     \\[3pt]   
B_{3,4} & =\ \ |\mathbf a \times \mathbf b |  +  \|\mathbf a  \pm  \mathbf b  \|^2     \\[3pt]   
          &=\ \ |\mathbf a \times \mathbf b| +\|\mathbf a\|^2 +\|\mathbf b\|^2  \pm 2 \mathbf a\bigcdot \mathbf b
\end{array}
\end{equation}
Consequently, the implied Apollonian disk packing is also integral. 
Moreover, for any vector
\begin{equation}
\label{eq:v}
\mathbf v = \alpha\mathbf a +\beta\mathbf b\,, \quad \gcd(\alpha,\beta)=1
\end{equation}
there is a disk in the major corona of the Apollonian packing with the curvature 
\begin{equation}
\label{eq:Bv}
B_{\mathbf v} = \ |\mathbf a \times \mathbf b |+ \|\mathbf v \|^2     
\end{equation}
If $\mathbf v$ and $\mathbf w$ are two such linear combinations and
$\mathbf v\times\mathbf w=\pm B \equiv |\mathbf a \times \mathbf b|$ then the disks 
$(B_0, B_{\mathbf v}, B_{\mathbf w})$ form a tricycle. 
Each of \eqref{eq:main} is a special case of \eqref{eq:Bv}.
}
\\
\\
{\bf Proof:} This follows by reversed reading of the equations \eqref{eq:BBB1} of
the previous proposition, supplemented with Theorems 2.2 and 2.5B: 
start with given spinors and build the Descartes quadruple.  
Note that setting the negative sign in $B_0=-|\mathbf a\times\mathbf b|$ assures the 
that disk $B_0$ is the major disk in a everted configuration. 
In other words, define $B =  |\mathbf a \times \mathbf b |  = |\/\det M\/| $,
and then
$$
(\;-B,\; B+ \|\mathbf a \|^2,\; B+ \|\mathbf b \|^2  ,\;   
B +  \|\mathbf a+ \mathbf b  \|^2\;)
$$
must satisfy the Descartes formula. 
Consult Figure \ref{fig:nice} for clarification.
Appendix A provides 
concrete coordinates for the disks, given the tangency spinors.
\qed

%
%

~

%
%
%
%

Now for some examples.
\\

\noindent
{\bf Example 1.}  Consider the following two spinors:
$$
\mathbf a = \begin{bmatrix}\;1\\\hbox{-2}\end{bmatrix}
\qquad
\mathbf b = \begin{bmatrix}2\\2\end{bmatrix}
\qquad
      M    = \begin{bmatrix}\;1\!&\!2\,\\\;\hbox{-}2\!&\!2\,\end{bmatrix}
$$
Then
$$
      \det M    = 6,  \quad 
      \|\mathbf a \|^2 = 1^2+(-2)^2 = 5,\quad
      \|\mathbf b\|^2 = 2^2+2^2 = 8,  \quad 
      \mathbf a \bigcdot\mathbf b = 
                                                    -2
$$
Hence the two Descartes configurations (for different sign in $\mathbf a\pm\mathbf b$) are
$$
(-6, 11,14, 15)
\qquad\hbox{and}\qquad
(-6, 11, 14, 23)
$$ 
as illustrated in Figure \ref{fig:nice}.
The figure shows also the Apollonian packing implied by the above Descartes configurations.

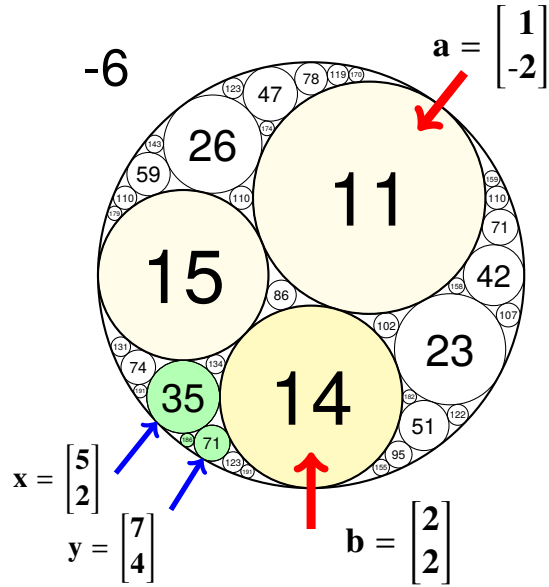
\begin{figure}[ht]
\centering
\begin{tikzpicture}[scale=17]

\draw [thick] (3/6,4/6) circle (1/6);
\draw [fill=yellow!10, thick] (6/11,8/11) circle (1/11);
\draw [fill=yellow!30, thick] (7/14,8/14) circle (1/14);
\draw [fill=yellow!10, thick] (6/15,10/15) circle (1/15);

\draw [fill=green!30, thin] (14/35,20/35) circle (1/35);
\draw [fill=green!30, thin] (30/71,38/71) circle (1/71);
\draw [fill=green!30, thin] (75/186,100/186) circle (1/186);
\foreach \a/\b/\c/\d/\s in {
 6/ 8/ 11/ 94/1.7,
 7/ 8/ 14/ 8/1.7,
 6/ 10/ 15/ 9/1.7,
14/ 14/ 23/ 17/1,
11/ 20/ 26/ 20/1,
14/ 20/ 35/ 17/.9,
27/ 28/ 42/ 36/.7,
22/ 38/ 47/ 41/.5,
30/ 28/ 51/ 33/.5,
22/ 44/ 59/ 41/.5
}
\draw (\a/\c,\b/\c) circle (1/\c)
        node [scale= 1.7*\s]  at (\a/\c,\b/\c)  {\sf \c\;}
;
\foreach \a/\b/\c/\d in {
30/ 38/ 71/ 33,
46/ 50/ 71/ 65,
27/ 44/ 74/ 36,
39/ 64/ 78/ 72,
41/ 56/ 86/ 56,
54/ 50/ 95/ 57,
57/ 64/ 102/ 72,
70/ 68/ 107/ 89,
49/ 80/ 110/ 80,
39/ 80/ 110/ 72,
71/ 80/ 110/ 104,
62/ 98/ 119/ 113,
75/ 68/ 122/ 84,
54/ 64/ 123/ 57,
54/ 100/ 123/ 105,
46/ 80/ 131/ 65,
57/ 80/ 134/ 72,
54/ 110/ 143/ 105,
86/ 80/ 155/ 89,
97/ 104/ 158/ 128,
102/ 118/ 159/ 153,
91/ 140/ 170/ 164,
81/ 136/ 174/ 144,
62/ 128/ 179/ 113,
105/ 104/ 182/ 120,
75/ 100/ 186/ 84,
86/ 98/ 191/ 89,
70/ 110/ 191/ 89
}
\draw (\a/\c,\b/\c) circle (1/\c)  
node [scale= 45/\c]  at (\a/\c,\b/\c)  {\sf \c \;}
;

\node [scale= 2]  at (.34, .83)  {\sf -6};

\draw [red,->, line width=1.2mm] (3/5+.012+.008, 4/5+.015+.01)--(3/5-.012-.004, 4/5-.015-.005);
\draw [red,->, line width=1.2mm] (4/8, 4/8-.032)--(4/8, 4/8+.028);

\draw node  [scale=1.2]  at (3/5+.042, 4/5+.045) 
                              {\color{black} $\mathbf a = \begin{bmatrix}\;\mathbf 1\\ \hbox{-}\mathbf 2\end{bmatrix}$};
\draw node  [scale=1.2]  at (4/8+.07, 4/8-.043) 
                              {\color{black} $\mathbf b = \begin{bmatrix}\mathbf 2 \\ \mathbf 2\end{bmatrix}$};

\draw [blue,->, line width=.7mm] (11/29-.032, 16/29-.04)--(11/29, 16/29);
\draw [blue,->, line width=.7mm] (27/65-.025, 34/65-.04)--(27/65, 34/65);
\draw node  [scale=1]  at (11/29-.077, 16/29-.045) 
                            {\color{black} $\mathbf x = \begin{bmatrix}\mathbf 5\\ \mathbf 2\end{bmatrix}$};
\draw node  [scale=1]  at (27/65-.07, 34/65-.07) 
                            {\color{black} $\mathbf y = \begin{bmatrix}\mathbf 7\\ \mathbf 4\end{bmatrix}$};

%
%
\end{tikzpicture}
\caption{Packing $(-6,11,14,15)$ and non-principal spinors}
\label{fig:notnice}
\end{figure}

~\\
{\bf Example 2}: Let us try another choice of   
integral vectors: 
$$
\mathbf x = \begin{bmatrix}5\\2\end{bmatrix}
\qquad
\mathbf y = \begin{bmatrix}7\\4\end{bmatrix}
\qquad
      M    = \begin{bmatrix}\,5\!&\!7\,\\\,2\!&\!4\,\end{bmatrix}
$$
Calculate
$$
      \det M    = 6,  \quad 
      \|\mathbf a \|^2 = 29,\quad
      \|\mathbf b\|^2 =  65,  \quad 
      \mathbf \|\mathbf a \pm\mathbf b\| = 43
$$
%
Hence the curvatures are determined
and the two Descartes configurations are
$$
(-6, 35, 164, 186)
\qquad\hbox{and}\qquad
(-6, 35, 164, 14)
$$ 
One may identify these disks in Figure \ref{fig:notnice} (green disks in the left lower region). 
Different Descartes quadruples may give the same Apollonian disk packing.
This example is to show that a random choice of integral spinors does not necessarily lead to 
a maximal Apollonian packing.
But their appropriate linear combination
will give the principal spinors and the corresponding maximal Descartes configuration:
$$
(\mathbf a,\, \mathbf b)  \ = \ (3\mathbf x - 2\mathbf y,\;  \mathbf y - \mathbf x)
$$

\newpage

\noindent
{\bf Example 3:}
Let us see what happens when one chooses two non-adjacent spinors in the Apollonian packing in Figure \ref{fig:nice}.
$$
\mathbf a = \begin{bmatrix}-1\\\;2\end{bmatrix}
\qquad
\mathbf x = \begin{bmatrix}5\\2\end{bmatrix}
$$
Since $\mathbf a\times \mathbf x = -12$ the spinors build a different 
Apollonian packing. 
The implied Descartes configuration is
$$
(-12,  17, 41,  46\pm 1) 
$$  
which will produce  a different Apollonian packing.


\subsection{Parametrizing formula interpreted}


The Descartes formula may be derived from the spinor theorems 2.1--2.5
(see  \cite{jk-spinor}).
One may start with the Pythagorean Theorem 
$\sin^2\varphi + \cos^2\varphi =1$, multiply it by absolute 
values of two vectors in $\mathbb R^2$ to get a known vector formula shown below.
Next, interpret the vectors as two adjacent spinors to get
the following correspondence:

{\large
\begin{equation}
\label{eq:big}
\boxed{
      \phantom{\Big|^{\big|}}\quad
        \underbrace{\phantom{\big|} (\mathbf a \times \mathbf b )^2}_{\mathlarger B^2} 
\ \ +\   \underbrace{\phantom{\big|}\,(\mathbf a \bigcdot \mathbf b )^2}_{\mathlarger \mu^2}  
\ \ \ = \ \ \   \underbrace{\phantom{\big|}\, \|\mathbf a\|^2\;\|\mathbf b\|^2}_{\mathlarger{ k\ \cdot \ n}\ }
\qquad\phantom{\Big|}
\atop \phantom{a}
}
\end{equation}
}

\noindent
This formula is valid for any tricycle.
But if one of the disks is assumed to be outer (i.e., of negative curvature), 
one  may readily recognize in the above equation the terms 
of the parametrizing formula of Theorem 1.1,  as shown in the above box.
Suddenly we have a spinor interpretation of the algebraic terms of the parametrizing formula.
Terms $k$ and $n$ are simply the norms squared of spinors and therefore they must be sums of squares.
Term $\mu$ turns out to be the inner product of spinors, and as such it represents the mid-circle $(B_0,B_1,B_2)$.
\\
\\
{\bf Remark:} To see that Equation \eqref{eq:big} is equivalent to Descartes' formula, 
replace the vector products with the expressions in terms of curvatures from Proposition 3.2 to \eqref{eq:big} to get 
\begin{equation}
\label{eq:BBBB}
    B_0^2 \ + \  \left(\frac{(B_0+B_3)-(B_0+B_1)-(B_0+B_2)}{2}\right)^2 \ = \ (B_0 + B_1)(B_0 + B_2)
\end{equation}
which reduces in a few steps to the known standard form \eqref{eq:Descartes}.
%
%
%


\subsection{Remark on Brahmaputra-Fibonacci identity}

From the number-theoretic point of view Equation \eqref{eq:big}
may be interpreted as an instance of Brahmaputra-Fibonacci identity:
$$
(a^2 + b^2)(c^2+d^2) = (ac+bd)^2 + (ad-bc)^2
$$
Thus the \eqref{eq:big} provides (1) a vector interpretation of the Brahmaputra-Fibonacci identity,
and (2) a geometric interpretation in terms of disks.
Moreover, it is clear that it is a form of Pythagorean Theorem!

\def\smalll{\scriptsize}
\begin{equation}
\begin{tikzpicture}[baseline=-0.8ex]
    \matrix (m) [ matrix of math nodes,
                         row sep=2em,
                         column sep=2.3em,
                         text height=3.8ex, text depth=3ex] 
   {
          {\hbox{\sf Pythagorean } \atop \hbox{theorem}}  \quad  
                  & \quad {\hbox{\sf Parametrizing} \atop \hbox{\sf equation}} \quad  
                  & \quad {\hbox{\sf Descartes}\atop\hbox{\sf formula}}  \\
           {~}  & \quad { \hbox{\sf Fibonacci-Brahmagupta}\atop \hbox{\sf Formula}} \qquad     \\
     };
    \path[<->]
        (m-1-1)   edge node[above] {$\hbox{ \sf\smalll  homo}$} (m-1-2)
        (m-1-2)   edge  (m-1-3)
        (m-2-2)   edge  (m-1-2)
        ;
\end{tikzpicture}   
\end{equation}

\noindent
{\bf Remark:} 
Note that the cross-product is defined only in dimensions 2, 3 and 7 
(the last two due to the existence of quaternions and octonions)
hence the Brahmaputra-Fibonacci identities have the corresponding generalized versions. 
\\

\noindent
{\bf Corollary:} 
{\sf 
The meaning of $\mu$ of the parametrizing formula is extended:
$$
\mu = \mathbf a \bigcdot \mathbf b  
       \ = \ \frac{1}{2} (B_3-B_2-B_1-B_0) 
       \ = \ \sqrt{B_0 B_1 + B_1 B_2 + B_2B_0}  =  T
$$        
where $T$ denotes curvature of the circle
through the points of tangency of circles 1,2, and 3.
}
\\

\noindent
{\bf Remark:} 
If the circle is to small (or, $\mu$ too big), the spinors $\mathbf a$ and $\mathbf b$ 
fail to be principal because the respective circles are too ``crumbled'' together, too small to like in Figure \ref{fig:notnice}.

\begin{figure}
{\large
$$
\boxed{
\begin{array}{cccccl}
~\\[-7pt]
 B_0^2  &+&   \frac{1}{4} \left(B_3\!-\!B_0\!-\!B_1\!-\!B_2\right)^2 &=& (B_0\! +\! B_1)(B_0 \!+\! B_2)
                                &\qquad\hbox{\normalsize\sf Descartes formula}   \\[7pt]
\mathlarger B^2   &+&   {\mathlarger \mu^2}  
                          & = &  {\mathlarger{ k\ \cdot \ n}\ }  
                          &\qquad\hbox{\normalsize \sf Parametrizing formula} \\[7pt]
(\mathbf a \times \mathbf b )^2 
                          &+&   (\mathbf a \bigcdot \mathbf b )^2 
                          &= &   \|\mathbf a\|^2\;\|\mathbf b\|^2
                          &\qquad\hbox{\normalsize\sf Spinorial formula} \\[7pt]
(ac+bd)^2 &+& (ad-bc)^2  &=  & (a^2 + b^2)(c^2+d^2) 
                          &\qquad\hbox{\normalsize\sf Fibonacci-Brahmaputra} \\[7pt]
\sin^2\varphi &+& \cos^2\varphi &=& 1             
                          &\qquad\hbox{\normalsize\sf Pythagorean Theorem}          \\
~
 \end{array}
}
$$
}
\caption{Five variations on a theme}
\label{fig:variations}
\end{figure}

\newpage

~

%


\newpage 
\section{Lattices}

Here we review some basic facts about lattices,
restricted to the case of two dimensions 
\cite{Lagrange,KM,Serre}. 
The two-dimensional integral lattice $\mathbb Z^2$ is called {\bf Diophantine plane}
or the {\bf standard integer lattice}. 

By an {\bf integral lattice} we will understand any 2-dimensional sublattice of $\mathbb Z^2$,
i.e., 
a $\mathbb Z$-module generated by two linearly independent vectors 
$\mathbf v$ and $\mathbf w$ of the Diophantine plane: 
\begin{equation}
\label{eq:lattice}
L(\mathbf v, \mathbf w) \ = \ \left\{\; \alpha \mathbf v+\beta \mathbf w\;|\; \alpha,\beta \in \mathbb Z\; \right\}
\end{equation}
Matrix whose columns are made by the basis vectors, $M = [\mathbf v\,\mathbf w]$,
allows one for abbreviated description of the lattice:
$$
L(\mathbf v,\mathbf w) \ = \ \{M\xi \;|\; \xi \in \mathbb Z^2\} \ = \  M\mathbb Z^2 
$$ 
The {\bf primitive cell} 
for a basis $\mathbf v$, $\mathbf w$ is the region 
$$
\{x\mathbf v+y\mathbf w \in \mathbb R^2\;|\; 0 \leq x,y \leq 1\;\}
$$
Placing a copy of the primitive cell at every lattice point produces a tiling of the plane.
Figure \ref{fig:twobases}  shows two choices of basis generating the same lattice 
and the corresponding two different primitive cells. 
The {\bf discriminant of a lattice} (called also the ``determinant'' of the lattice)
is defined as the area of the primitive cell
$$
\det(L) \ = \ 
         |\, \mathbf v\times \mathbf w\,| = |\, \det M\,|
$$

 It does not depend on  the choice of basis.
Two bases define the same lattice if they are the  image one of the other via unimodular matrix, 
i.e. transformation from the group 
$$
T\in \hbox{SL}(2,\mathbb Z)
$$

\def\nmaly{0.37}  
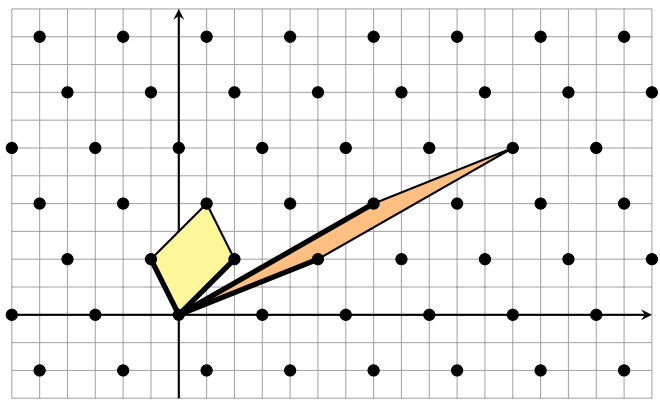
\begin{figure}[H]
\centering
\begin{tikzpicture}[scale=\nmaly, rotate=0, shift={(0,0cm)}]  
\clip (-6.5,-4) rectangle (17.5,12);
\def\rr{.2}
\def\Ax{5}; \def\Ay{2};
\def\Bx{7}; \def\By{4};

\def\Cx{2}; \def\Cy{2};
\def\Dx{-1}; \def\Dy{2};

\draw[thin, color=gray!70] (-6,-3) grid (17,11);
\draw[-stealth,  thick] (-6,0) -- (17,0);   
\draw[-stealth,  thick] (0, -3) -- (0, 11);   

\draw [thick, fill=orange!50, fill opacity=.7]  (0,0) -- (\Ax,\Ay) --  (\Ax+\Bx, \Ay + \By)  -- (\Bx,\By) -- cycle;
\draw [thick, fill=yellow!50, fill opacity=.7]  (0,0) -- (\Cx,\Cy) --  (\Cx+\Dx, \Cy + \Dy)  -- (\Dx,\Dy) -- cycle;

\draw [fill=white] (-8,-2) circle  (\rr);
\draw [fill=black] (-5,-2) circle  (\rr);
\draw [fill=black] (-2,-2) circle  (\rr);
\draw [fill=black] (1,-2) circle  (\rr);
\draw [fill=black] (4,-2) circle  (\rr);
\draw [fill=black] (7,-2) circle  (\rr);
\draw [fill=black] (10,-2) circle  (\rr);
\draw [fill=black] (13,-2) circle  (\rr);
\draw [fill=black] (16,-2) circle  (\rr);

\draw [fill=black] (-6,0) circle  (\rr);
\draw [fill=black] (-3,0) circle  (\rr);
\draw [fill=black] (0,0) circle  (\rr);
\draw [fill=black] (3,0) circle  (\rr);
\draw [fill=black] (6,0) circle  (\rr);
\draw [fill=black] (9,0) circle  (\rr);
\draw [fill=black] (12,0) circle  (\rr);
\draw [fill=black] (15,0) circle  (\rr);

\draw [fill=black] (-7,2) circle  (\rr);
\draw [fill=black] (-4,2) circle  (\rr);
\draw [fill=black] (-1,2) circle  (\rr);
\draw [fill=black] (2,2) circle  (\rr);
\draw [fill=black] (5,2) circle  (\rr);
\draw [fill=black] (8,2) circle  (\rr);
\draw [fill=black] (11,2) circle  (\rr);
\draw [fill=black] (14,2) circle  (\rr);
\draw [fill=black] (17,2) circle  (\rr);

\draw [fill=black] (-8,4) circle  (\rr);
\draw [fill=black] (-5,4) circle  (\rr);
\draw [fill=black] (-2,4) circle  (\rr);
\draw [fill=black] (1,4) circle  (\rr);
\draw [fill=black] (4,4) circle  (\rr);
\draw [fill=black] (7,4) circle  (\rr);
\draw [fill=black] (10,4) circle  (\rr);
\draw [fill=black] (13,4) circle  (\rr);
\draw [fill=black] (16,4) circle  (\rr);

\draw [fill=black] (-6,6) circle  (\rr);
\draw [fill=black] (-3,6) circle  (\rr);
\draw [fill=black] (0,6) circle  (\rr);
\draw [fill=black] (3,6) circle  (\rr);
\draw [fill=black] (6,6) circle  (\rr);
\draw [fill=black] (9,6) circle  (\rr);
\draw [fill=black] (12,6) circle  (\rr);
\draw [fill=black] (15,6) circle  (\rr);
\draw [fill=black] (18,6) circle  (\rr);

\draw [fill=black] (-7,8) circle  (\rr);
\draw [fill=black] (-4,8) circle  (\rr);
\draw [fill=black] (-1,8) circle  (\rr);
\draw [fill=black] (2,8) circle  (\rr);
\draw [fill=black] (5,8) circle  (\rr);
\draw [fill=black] (8,8) circle  (\rr);
\draw [fill=black] (11,8) circle  (\rr);
\draw [fill=black] (14,8) circle  (\rr);
\draw [fill=black] (17,8) circle  (\rr);

\draw [fill=black] (-5,10) circle  (\rr);
\draw [fill=black] (-2,10) circle  (\rr);
\draw [fill=black] (1,10) circle  (\rr);
\draw [fill=black] (4,10) circle  (\rr);
\draw [fill=black] (7,10) circle  (\rr);
\draw [fill=black] (10,10) circle  (\rr);
\draw [fill=black] (13,10) circle  (\rr);
\draw [fill=black] (16,10) circle  (\rr);
\draw [fill=black] (18,10) circle  (\rr);


\draw [line width=2pt] (0,0) -- (\Ax,\Ay);
\draw [line width=2pt] (0,0) -- (\Bx,\By);
\draw [line width=2pt] (0,0) -- (\Cx,\Cy);
\draw [line width=2pt] (0,0) -- (\Dx,\Dy);

\end{tikzpicture}
\caption{Two bases in the same lattice}
\label{fig:twobases}
\end{figure}

A basis is called {\bf principal} if it consists of the shortest vectors of the lattice.
One of the bases in Figure \ref{fig:twobases} is principal.
Note that such a reduced basis may be always chosen so that 
the basis vectors make a right or acute angle:
$$
\mathbf b_1 \bigcdot \mathbf b_2 \geq 0
$$
Indeed, you may always replace one of the basis vectors by its negative, as it has the same length,
tochange the sign of the above inner produc.
We shall include this to the definition; A basis is {\bf principal} if it is ordered, 
reduced and makes an acute angle: 
$$
\|\mathbf b_1\|
\leq
\|\mathbf b_2\| 
\leq 
\|\mathbf  b_1\pm \mathbf  b_2\| 
\qquad\hbox{and}\qquad 
\mathbf b_1\bigcdot \mathbf b_2 \geq 1
$$

To find a principal basis for a given lattice, one performs a {\bf basis reduction process}.
Start with any basis $\mathbf v,\, \mathbf w$
and replace the longer vector, say $\mathbf w$, by $\mathbf w' = \mathbf w\pm \mathbf v$ to 
get a new basis $(\mathbf w',\, \mathbf v)$.  
Repeat until such a move does not reduce the size of vectors any more.
This is an analogue of the Euler's algorithm for computing the greatest common divisor. 
An obvious (in 2D case) characterization follows: 
a basis $\mathbf b_1, \mathbf b_2$ of a lattice is Lagrange-Gauss reduced 
if
$
\|\mathbf b_2+q\mathbf b_1\| \geq \|\mathbf b_2\| \geq \|\mathbf b_1\|
$
for every $q\in \mathbb Z$. 
In the case of 2-dimensional basis this coincides with the principal basis.
\\

We shall call two integral lattices {\bf similar} if one may be obtained from the other by 
any composition of rotation, reflection and dilation.  Figure \ref{fig:similar} shows an example of three 
mutually similar lattices. (The numbers denote the discriminant of the lattices.)
\\

\def\nmaly{.37}
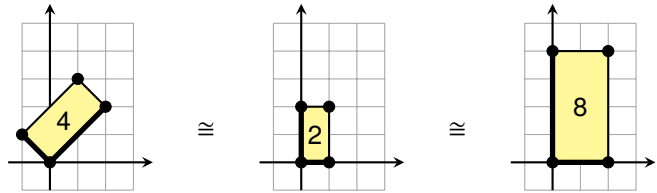
\begin{figure}[H]
\centering
\newcommand{\romb}[7]{%
\begin{tikzpicture}[scale=\nmaly, rotate=0, shift={(0,0cm)}]  
\draw[thin, color=gray!70] (-1,-1) grid (#1,#2);
\draw[-stealth,  thick] (-1.5,0) -- (#1 +.7,0);   
\draw[-stealth,  thick] (0, -1) -- (0, #2+.7);   
\draw [thick, fill=yellow!50, fill opacity=.7]  (0,0) -- (#3,#4) --  (#3+#5, #4+#6)  -- (#5,#6) -- cycle;
\draw [fill=black] (0,0) circle  (.2);
\draw [fill=black]  (#3,#4) circle  (.2);
\draw [fill=black]  (#5,#6) circle  (.2);
\draw [fill=black]  (#3+#5, #4+#6) circle  (.2);
\draw [line width=2pt] (0,0) -- (#3,#4);
\draw [line width=2pt] (0,0) -- (#5,#6);
\draw node at  (.5*#3+.5*#5,.5*#4+.5*#6) [scale=1]{\sf #7};
\end{tikzpicture}%
}
$$
\romb{3}{5}{-1}{1}{2}{2}{4}
\quad\begin{array}{c}\cong\\~\\~\\~\end{array}\quad
\romb{3}{5}{1}{0}{0}{2}{2}
\quad\begin{array}{c}\cong\\~\\~\\~\end{array}\quad
\romb{3}{5}{2}{0}{0}{4}{8}
$$
\vspace{-.5in}
\caption{Similar lattices}
\label{fig:similar}
\end{figure}

%
%
%

For every integer $n\in \mathbb N$ there exist a finite number of integral lattices of discriminant $n$.
Figure \ref{fig:lattices} shows them for $n= 1,..., 9$, and $n=12$.

~

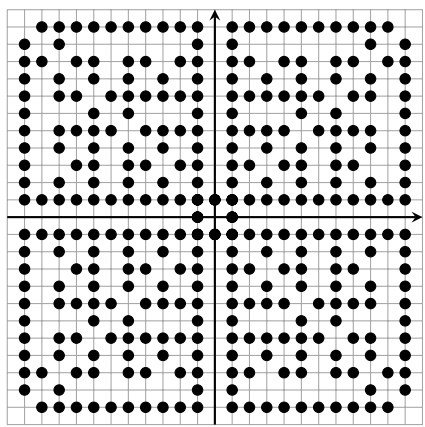
\begin{figure}[H]
\centering
\begin{tikzpicture}[scale=.23, rotate=0, shift={(0,0cm)}]  
\def\N{12}; \def\rr{.3}
\clip (-\N-.03,-\N-.03) rectangle (\N+.03,\N+.03);
\draw[thin, color=gray!70] (-\N,-\N) grid (\N,\N);
\draw[-stealth,  thick] (-\N,0) -- (\N,0);   
\draw[-stealth,  thick] (0, -\N) -- (0, \N);   
\foreach \a/\b in {
1 / 0, 1/1,
2/1, 
3/1,3/2, 
4/1, 4/3, 
5/1, 5/2, 5/3, 5/4,
6/1, 6/5, 
7/1, 7/2, 7/3 , 7/4, 7/5, 7/6,
8/1, 8/3, 8/5, 8/7 ,
9/1, 9/2,9/4,9/5,9/7,9/8,
10/1,10/310/7,10/9,
11/1,11/2,11/3,11/4,11/5,11/6,11/7,11/8,11/9,11/10
}
\draw [fill=black]	(\a,\b) circle  (\rr)
				(\a,-\b) circle  (\rr)
 				(-\a,\b) circle  (\rr)
 				(-\a,-\b) circle  (\rr)
				(\b,\a) circle  (\rr)
  				(\b,-\a) circle  (\rr)
 				(-\b,\a) circle  (\rr)
 				(-\b,-\a) circle  (\rr)
;
\end{tikzpicture}
\caption{Euclidean mosaic (a subset of the standard lattice)}
\label{fig:u-lattice}
\end{figure}

We shall introduce yet another term: 
The 2-dimensional {\bf Euclidean mosaic} is the set of vectors (points) in the Diophantine plane
whose components are coprimes:
$$
\mathbb Z^2_{\circ} = \{ [n,k]^T \;|\; \gcd(n,k)=1\;\}
$$
Equivalently, it is the orbit of $[1,0]$ via the action of the unimodular group $\SL^{\pm}(2,\mathbb Z)$,
the group of $2\!\times\! 2$ matrices of determinant equal to $\pm 1$.
Euclidean mosaic is {\it not} a lattice, as it is not a group
(see Figure  \ref{fig:u-lattice}). 
{\bf Mosaic of a lattice} $L(\mathbf v,\mathbf w)$ 
is defined analogously as the subset
 $$
 L_{\circ}(\mathbf v,\mathbf w)=\{\alpha \mathbf v +\beta\mathbf w \;|\; \gcd(\alpha,\beta)=1\}
 \  = \ 
 \{ \, M\, \mathbf \xi \;|\; \mathbf \xi \in \mathbb Z^2_o\;\}
 $$
The reason for introducing this term is the fact that due to Theorem 2.5B, 
the spinors of a given corona form a mosaic in the lattice defined by spinors.

\def\nmaly{0.37}  
\begin{figure}[H]
\centering
\newcommand{\romb}[7]{%
\begin{tikzpicture}[scale=\nmaly, rotate=0, shift={(0,0cm)}]  
\draw[thin, color=gray!70] (-1,-1) grid (#1,#2);
\draw[-stealth,  thick] (-1.5,0) -- (#1 +.7,0);   
\draw[-stealth,  thick] (0, -1) -- (0, #2+.7);   
\draw [thick, fill=yellow!50, fill opacity=.7]  (0,0) -- (#3,#4) --  (#3+#5, #4+#6)  -- (#5,#6) -- cycle;
\draw [fill=black] (0,0) circle  (.2);
\draw [fill=black]  (#3,#4) circle  (.2);
\draw [fill=black]  (#5,#6) circle  (.2);
\draw [fill=black]  (#3+#5, #4+#6) circle  (.2);
\draw [line width=2pt] (0,0) -- (#3,#4);
\draw [line width=2pt] (0,0) -- (#5,#6);
\draw node at  (.5*#3+.5*#5,.5*#4+.5*#6) [scale=1]{\sf #7};
\end{tikzpicture}%
}

\romb{2}{4}{1}{0}{0}{1}{1}
\romb{2}{4}{1}{0}{0}{2}{2}
\romb{2}{4}{1}{0}{0}{3}{3}
\romb{2}{4}{1}{1}{-1}{2}{3}
\romb{2}{4}{1}{0}{0}{4}{4}
\romb{3}{4}{2}{0}{1}{2}{4}

~

\romb{3}{6}{1}{0}{0}{5}{5}
\romb{3}{6}{2}{3}{-1}{1}{5}
\romb{3}{6}{1}{0}{0}{6}{6}
\romb{3}{6}{2}{0}{0}{3}{6}
\romb{3}{6}{2}{2}{-1}{2}{6}

~

\romb{3}{8}{1}{0}{0}{7}{7}
\romb{3}{8}{3}{4}{-1}{1}{7}
\romb{3}{8}{2}{1}{-1}{3}{7}
\romb{3}{8}{1}{0}{0}{8}{8}
\romb{3}{8}{3}{2}{-1}{2}{8}
\romb{3}{8}{2}{0}{1}{4}{8}

~

\romb{3}{9}{1}{0}{0}{9}{9}
\romb{3}{9}{4}{5}{-1}{1}{9}
\romb{3}{9}{2}{1}{-1}{4}{9}
\romb{3}{9}{3}{1}{0}{3}{9}

\caption{The first 21 irreducible lattices (and one more)}
\label{fig:lattices}
\end{figure}
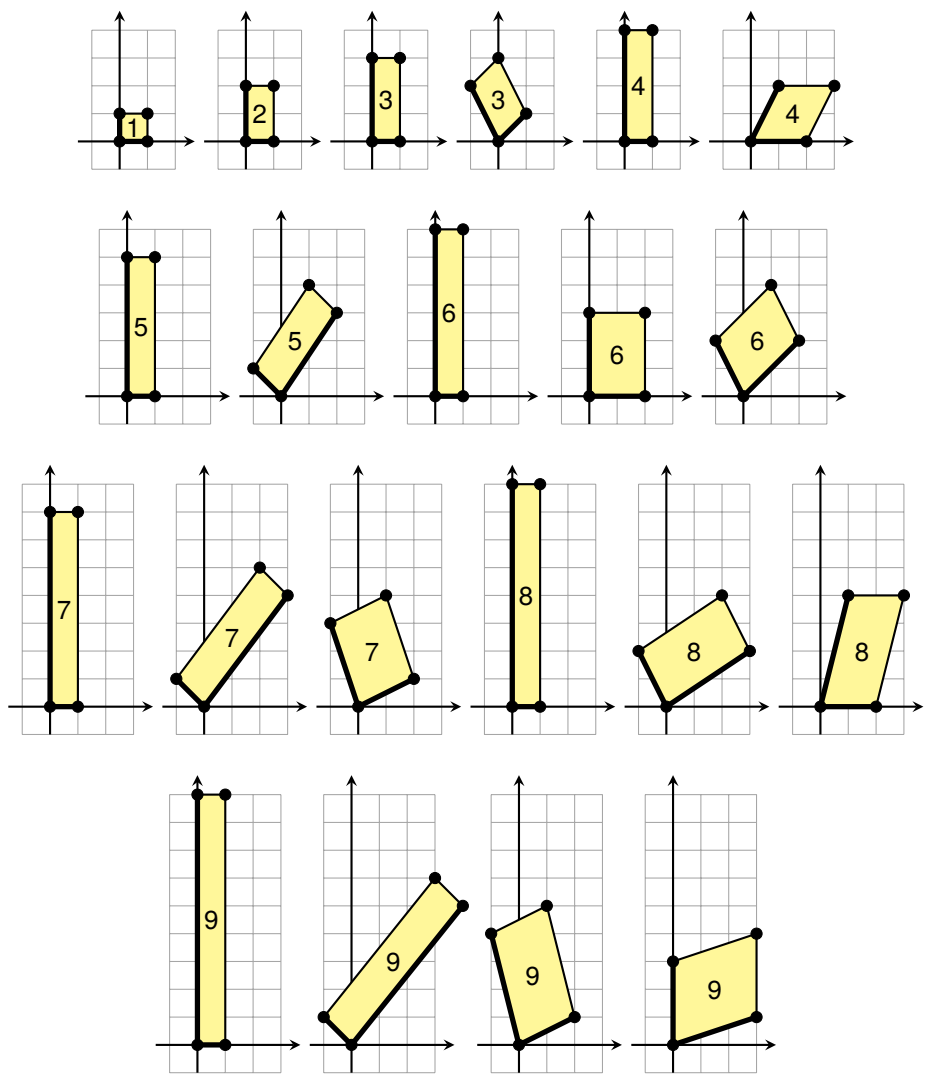

%

\subsection{Lattices and Apollonian packings}

With the language of lattices, we may now restate the content of Proposition 3.1 and 3.2 as follows:

~\\
{\bf Proposition 4.1.}
{\sf
The  corona spinors  in an Apollonian packing 
form a mosaic in the lattice defined by any pair 
of the adjacent spinors in the corona.
}
\\
\\
{\bf Proof:}  Any two adjacent spinors $\mathbf v, \mathbf w$ in the corona 
determine a lattice $L(\mathbf v,\mathbf w)$.
By Theorem 2.5B, any other spinor in the corona is a linear  combination of the two, 
hence is an elements of the same lattice.
But since the combinations involve coprime coefficients, they will make a subset of the lattice,
coinciding with its Euclidean mosaic.
\qed
\\

The very existence of the lattice allows one to find its principal basis -- two shortest vectors in the lattice -- 
through the reduction process described the previous subsection.
They point to (or from) the greatest circles in the corona, and consequently
indicate the maximal Descartes configuration in the Apollonian packing.
Thus we get a connection with the parametrizing formula of Theorem 1.1,
together with the interpretation of $k$, $n$, and $\mu$ given in \eqref{eq:big},
now applied to the principal spinors.
\\

This establishes also uniqueness of the classification
via the lattices up to similarity.
Indeed, rotations and reflections do not change the values   
of $k$, $n$ and $\mu$, while scaling would 
lead to quadruples of curvatures with a common divisor.
\\

In conclusion, the list of lattices, as illustrated in Figure \ref{fig:lattices}
is equivalent to the classification of the irreducible Apollonian packings. 
\\

\noindent
{\bf Remark 1:}
To achieve uniqueness of the spinor representation of lattices,
we may demand an extended list of conditions: 
\\[-10pt]

(1) $\|\mathbf a\|^2 \geq \|\mathbf b\|^2$ \ --- vector $\mathbf a$ is shorter than $\mathbf b$.

(2) $0 < a_2 < a_1$  \ --- vector $\mathbf a$ is in the first quadrant, close to the x-axis.

(3) $b_2 \geq 0$ \ --- vector $\mathbf b$ is in the first or second quadrant.

(4) $\mathbf a\bigcdot \mathbf b\geq0$ \ --- the vectors form an acute or right angle.
\\

\noindent
{\bf Remark 2:}
The above is true for any corona on any disk.
Only when our attention was on classification of the integral Apollonian packing, 
we concentrate on the outer disk with negative curvature. 
\\

\noindent
{\bf Proposition 3:}
{\sf
Every tangency spinor $\mathbf u$ in the  Apollonian packing is a linear combination
$$
\mathbf u \ = \ \alpha \mathbf v +\beta \mathbf w +\gamma \mathbf v^{\symp} +\gamma \mathbf w^{\symp}
$$
where $\alpha,\beta,\gamma,\delta \in \mathbb Z$, 
and $\mathbf v$ and $\mathbf w$ is a pair of adjacent spinors at {\it any} disk in the packing.
In particular, if $\mathbf v$ and $\mathbf w$ are integral, so are all spinors in the packing.
}
\\

\noindent
{\bf Proof:}  This is a simple consequence of Theorems 2.4 and 2.5, which allow one 
to extend the values of spinors throughout the packing.
\QED

\newpage

\subsection{Derivation of the constraints}

The conditions on $k$, $n$ and $\mu$ for the parametrizing formula in Theorem 1.1 
obtain a simple geometric justification when considered in the lattice context.
Recall that two spinors that determine a lattice may not correspond to the maximal Descartes configuration.
But they may be adjusted by the reduction process, i.e., a sequence of mutual subtraction / additions,  
until none can be replaced by a shorter one.

~\\
{\bf Proposition 4.1:}
{\sf 
Principality of spinors defining a Descartes configuration implies the conditions \eqref{eq:conditions} 
of Theorem \ref{thm:1}.
}
\\

\noindent
{\bf Proof:}
Let us assume that $\mathbf v$ and $\mathbf w$ are two spinors forming a principal basis of the lattice 
defining the corona of an Apollonian packing.
Recall  notation $k=\|\mathbf v\|^2$ and $n=\|\mathbf w\|^2$, and assume
the order
\begin{equation}
\label{eq:condition0}
\boxed{
 ~ |\mathbf v|\ \leq \ |\mathbf w| 
 \phantom{\big|} \; }
\end{equation}
(or $k\leq n$).  
Principality implies that the longer vector $\mathbf w$ cannot be replaced by 
$\mathbf w - \mathbf v$ to make the basis vectors shorter.
In other words: 
$$
|\mathbf w-\mathbf v| \geq |\mathbf w|
$$
Squaring both sides,
$$  
|\mathbf w|^2 - 2\mathbf w\bigcdot \mathbf v + |\mathbf v|^2 \ \geq \  |\mathbf w|^2
$$
leads to
\begin{equation}
\label{eq:condition1}
\boxed{ ~ 
2\mathbf w\bigcdot \mathbf v  \leq  |\mathbf v|^2
\phantom{\big|}\; }
\end{equation}
Now, this and the assumption $|\mathbf v| \leq |\mathbf w|$ imply 
$$
2\,|\mathbf w|\!\cdot\! |\mathbf v|\cos\varphi \ \leq \ |\mathbf v|^2 \ \leq \ |\mathbf v|\!\cdot\! |\mathbf w|\,,
$$
and, after cancellation,  
$$
\cos\varphi \ \leq \ \frac{1}{2}\,.
$$
In other words, $\varphi \geq 60^{\circ}$, or 
$$
\tan\varphi \ \geq \ \sqrt{3}\,.
$$
This implies $\sin\varphi \geq \sqrt{3}\cos\varphi$.
Multiplying both sides by $|\mathbf w|\!\cdot\! |\mathbf v|$ we get
\begin{equation}
\label{eq:condition2}
\boxed{ ~
\mathbf v\times \mathbf w \ \geq \ \sqrt{3}\, \mathbf v\bigcdot \mathbf w
\phantom{\big|}\; }
\end{equation}
In terms of notation of \eqref{eq:big}:
$$
k=\|\mathbf v\|^2, \qquad n=\|\mathbf w\|^2, 
\qquad B=|\mathbf v\times\mathbf w|, \quad \mu = \mathbf v\cdot \mathbf w
$$
the above properties \eqref{eq:condition0}, \eqref{eq:condition1}, and  \eqref{eq:condition2},  translate into: 
$$
2\mu \leq k \leq n \qquad \hbox{and} \qquad B/\sqrt{3} \ \geq \ \mu \,,
$$
as stated.
\qed

~

\noindent
{\bf Remark:}
Intuitively, condition $\mu<B/\sqrt{3}$ 
is equivalent to the geometric requirement that the mid-circle $(B_0,B_1,B_2)$ be not too small; 
otherwise, a circle greater than any of the two, $B_1$ or $B_2$ will be possible in the corona
of $B_0$.

\subsection{Curious examples of Apollonian strip and of unbounded ``golden'' packing}

A few examples that work, even if they seem degenerated.

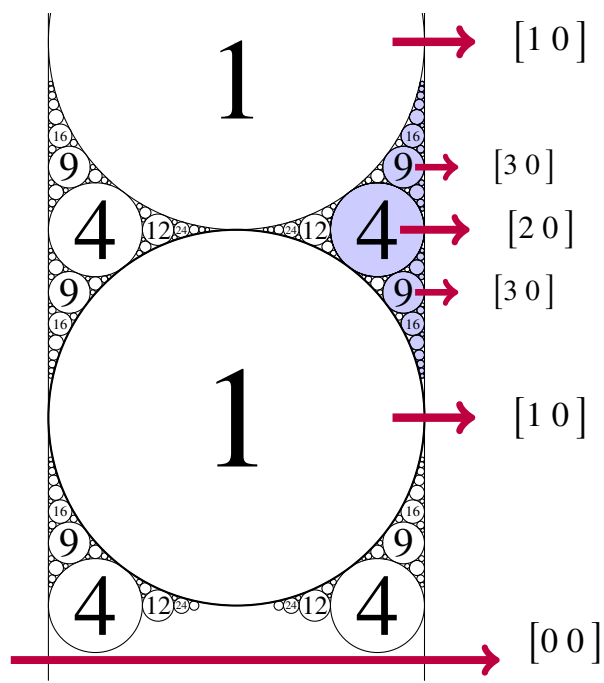
\begin{figure}[H]
\centering
\begin{tikzpicture}[scale=2.5]  

\clip (-1.4,-1.4) rectangle (2,2.15);
\foreach \a/\b/\c/\d in {
3/  4/ 4 ,
8/    6/   9 ,
24/  20/  25,
24/  10/  25,
15/   8/   16, 
35/	12/	36,
48/	42/	49,
48/	28/	49,
48/	14/	49,
63/	48/	64,	
63/	16/	64,	
80/	72/	81,	
80/	36/	81,	
80/	18/	81	
}
\draw  [fill=blue!20] (\a/\c,\b/\c) circle (1/\c)    
                                            (\a/\c,2-\b/\c) circle (1/\c)   
;
\draw (1,-1.4) -- (1,3);
\draw (-1,-1.4) -- (-1,3);
\draw (0,0) circle (1);
\draw [thick] (0,0) circle (1)
          node  [scale=5]  at (0,0) {$1$};
\draw (0,2) circle (1)
          node  [scale=4]  at (0,1.8) {$1$};
\foreach \a/\b/\c/\s in {
3/  4/ 4/    3,
5/ 12/12/  .9, 
7/ 24/24/ .4
}
\draw (\a/\c,\b/\c) circle (1/\c)    (\a/\c,-\b/\c) circle (1/\c)
          (-\a/\c,\b/\c) circle (1/\c)    (-\a/\c,-\b/\c) circle (1/\c)
        node [scale= \s]  at (\a/\c,\b/\c)  {$\c$}
        node [scale= \s]  at (\a/\c,-\b/\c) {$\c$}
        node [scale= \s]  at (-\a/\c,\b/\c) {$\c$}
        node [scale= \s]  at (-\a/\c,-\b/\c) {$\c$};
\foreach \a/\b/\c/\s in {
9/ 40/40/ 0
}
\draw (\a/\c,\b/\c) circle (1/\c)    (\a/\c,-\b/\c) circle (1/\c)
          (-\a/\c,\b/\c) circle (1/\c)    (-\a/\c,-\b/\c) circle (1/\c)
;
\foreach \a/\b/\c in {
11/	60	/60,
13/	84/	84
}
\draw (\a/\c,\b/\c) circle (1/\c)  
           (-\a/\c,\b/\c) circle (1/\c) 
;

\foreach \a/\b/\c/\d in {
8/    6/   9     / 1.4,
15/   8/   16  / .5
}
\draw (\a/\c,\b/\c) circle (1/\c)       (-\a/\c,\b/\c) circle (1/\c)
          (\a/\c,-\b/\c) circle (1/\c)       (-\a/\c,-\b/\c) circle (1/\c)
          (\a/\c,2-\b/\c) circle (1/\c)       (-\a/\c,2-\b/\c) circle (1/\c)
node [scale=\d]  at (\a/\c,\b/\c)    {$\c$}
node [scale=\d]  at (-\a/\c,\b/\c)   {$\c$}
node [scale=\d]  at (\a/\c,-\b/\c)   {$\c$}
node [scale=\d]  at (-\a/\c,-\b/\c)  {$\c$}
node [scale=\d]  at (\a/\c,2-\b/\c)    {$\c$}
node [scale=\d]  at (-\a/\c,2-\b/\c)   {$\c$};

\foreach \a /  \b / \c   in {
24/20/25,   24/10/25,    21/20/28,
16/30/33,   35/12/36,    48/42/49,
48/28/49,   48/14/49,    45/28/52
}
\draw (\a/\c,\b/\c) circle (1/\c)          (-\a/\c,\b/\c) circle (1/\c)
          (\a/\c,-\b/\c) circle (1/\c)         (-\a/\c,-\b/\c) circle (1/\c)
          (\a/\c,2-\b/\c) circle (1/\c)          (-\a/\c,2-\b/\c) circle (1/\c)
;


\foreach \a /  \b / \c   in {
40/	42/	57,	33/	56/	64,	63/	48/	64,	63/	16/	64,	55/	48/	72,	24/	70/	73,
69/	60/	76,	80/	72/	81,	64/	60/	81,	80/	36/	81,	80/	18/	81,	77/	36/	84
39/	80/	88,	65/	72/	96,	48/	92/	97,	99/	60/	100
}
\draw (\a/\c,\b/\c) circle (1/\c)          (-\a/\c,\b/\c) circle (1/\c)
          (\a/\c,-\b/\c) circle (1/\c)         (-\a/\c,-\b/\c) circle (1/\c)
          (\a/\c,2-\b/\c) circle (1/\c)          (-\a/\c,2-\b/\c) circle (1/\c);

\draw [->, line width=3.2pt, purple]  (.83,0) -- (1.27,0);
\node at (1.67, 0) [scale=1.2] {$\begin{bmatrix} 1 \  0\, \end{bmatrix}$};

\draw [->, line width=3.2pt, purple]  (.83,2) -- (1.27,2);
\node at (1.67, 2) [scale=1.2] {$\begin{bmatrix} 1 \ 0\, \end{bmatrix}$};

\draw [->, line width=3pt, purple] (.87,1) -- (1.24,1);
\node at (1.62, 1) [scale=1.1] {$\begin{bmatrix} 2\ 0\, \end{bmatrix}$};

\draw [->, line width=2.5pt, purple]  (.95,.667) -- (1.18,.667);
\node at (1.53, .667) {$\begin{bmatrix} 3\  0\, \end{bmatrix}$};

\draw [->, line width=2.5pt, purple]  (.95, 1.333) -- (1.18, 1.333);
\node at (1.53, 1.333) {$\begin{bmatrix} 3 \  0\, \end{bmatrix}$};

\draw [->, line width=3pt, purple]  (-1.2, -1.29) -- (1.4, -1.29);
\node at (1.75, -1.2) [scale=1.2] {$\begin{bmatrix} 0 \ 0\, \end{bmatrix}$};

\end{tikzpicture}
\caption{Ford circles as a corona in the Apollonian Belt}
\label{fig:strip}
\end{figure}

~\\
{\bf Example 4:}
Interestingly, the initial spinors do not need to be linearly independent
to determine an Apollonian packing.  
Start with 
$$
\mathbf a = \begin{bmatrix}1\\0\end{bmatrix}\,,
\qquad
\mathbf b = \begin{bmatrix}1\\0\end{bmatrix} \,,
\qquad
      M    = \begin{bmatrix}\,1\!&\!1\,\\\,0\!&\!0\,\end{bmatrix} \,.
$$
Calculate
$\|\mathbf a \|^2 = 1$, 
$\|\mathbf b \|^2 = 1$, 
$|\mathbf a \times \mathbf b| = $, 
and
$\mathbf a\bigcdot \mathbf b  = 1$. 
The two complementary Descartes configurations are therefore 
$$
(0,  1, 1,4) \qquad \hbox{and}\qquad (0,1,1,0)
$$  
One readily recognizes that they generate the Apollonian strip,
see Figure \ref{fig:strip}.
The integral sublattice is here one-dimensional, $\mathbb Z\subset \mathbb Z^2$.

~\\
{\bf Example 5:}
One does not even need a requirement of non-zero vectors!
For the sake of an experiment, let us try  
$$
\mathbf a = \begin{bmatrix}0\\0\end{bmatrix} \,,
\qquad
\mathbf b = \begin{bmatrix}1\\0\end{bmatrix} \,,
\qquad
      M    = \begin{bmatrix}\,0\!&\!1\,\\   
                                       \,0\!&\!0\,\end{bmatrix}   \,.
$$
Then 
$\|\mathbf a \|^2 = 0$, 
$\|\mathbf b \|^2 = 1$, 
$|\mathbf a \times \mathbf b| =0 $, 
and
$\mathbf a\bigcdot \mathbf b  = 0$,  
and the two Descartes configurations are 
$$
(0,  0, 1, 1) \qquad \hbox{and}\qquad (0,0,1,1)
$$  
which, as in the previous example,  may be completed to the Apollonian strip.
One should consider this case as the reduced, principal, basis of Example 4.
Tangent spinor $[0, 0]^T$ connects the straight lines (the half-spaces),
considered as tangent at infinity, see Figure~\ref{fig:strip}.
Clearly, this  one-dimensional lattice needs to be included into 
the list of lattices corresponding to the integral Apollonian packings.

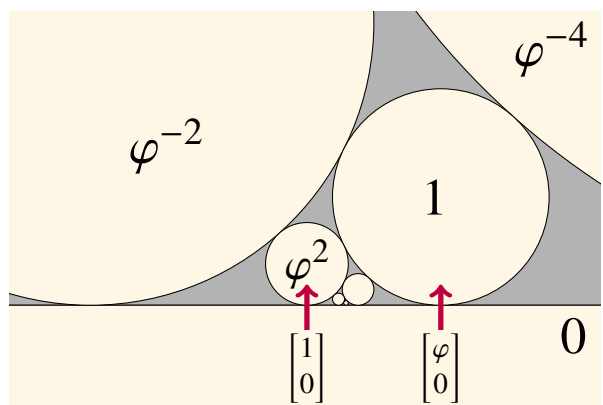
\begin{figure}[H]
\begin{center}
\begin{tikzpicture}[scale=1.1]

\clip (-3.6, -1.28) rectangle  (3.55,3.55);

\draw [fill=black!30] (-4,-2) rectangle (4,4);
\draw (-4,0) -- (4,0);

\draw [fill=gold!10] (-4,-2) rectangle (4, 0);

\foreach \a/\b   in {
8.472135955/	8.972135955,
-2.618033989/	3.427050983,
1.618033989/	1.309016994,
0/	0.5,
0.618033989/	0.190983006,
0.381966011/	0.072949017,
0.472135955/	0.027864045
}
\draw[fill=gold!10]  (\a,\b) circle (\b)
;
\node at (2.8,  3.1) [align=left, scale=1.8] {\ \ $\varphi^{-4}$}; 
\node at (-2.618033989+1,	3.427050983-1.5) [align=left, scale=1.8] {$\varphi^{-2}$ }; 
\node at (1.618033989, 1.309016994)  [align=left, scale=1.8] {$1$ }; 
\node at (0,	0.5) [align=left, scale=1.6] {$\; \varphi^{2}$ }; 

\draw [->, line width=2pt, purple]  (0, -.3) -- (0, .25);
\node at (0, -.76) [scale=1] {$\begin{bmatrix} 1 \\ 0 \end{bmatrix}$};

\draw [->, line width=2pt, purple]  (1.618, -.3) -- (1.618033989, .25);
\node at (1.618033989, -.76) [scale=1] {$\begin{bmatrix} \varphi \\ 0 \end{bmatrix}$};

\node at (3.22, -.37) [scale=1.8] {$0$};

\end{tikzpicture}
\end{center}
\caption{Upper half-plane Apollonian disk packing.  
To get actual orientation for the presented spinors, 
rotate the image by 90 degrees.}
\label{fig:halfplane}
\end{figure}

~\\
{\bf Example 6:}
Start with 
$$
\mathbf a = \begin{bmatrix}1\\0\end{bmatrix}
\qquad
\mathbf b = \begin{bmatrix}\varphi\\0\end{bmatrix}
\qquad
      M    = \begin{bmatrix}\,1\!&\!\varphi\,\\\,0\!&\!0\,\end{bmatrix}
$$
where $\varphi=(1+\sqrt{5})/2$ is the golden ratio.
Then the two Descartes configurations are 
$$
(0, \, 1,\,  \varphi^2,\, \varphi^4) \qquad \hbox{and}\qquad (0,\, 1,\, \varphi^2,\, \varphi^{-2})
$$  
As shown in \cite{jk-gold}, this corresponds to an unbounded Apollonian packing 
which does not contain a disk of  minimum positive curvature
(see Figure~\ref{fig:halfplane}).  
Note that that the process of reduction would be a never-ending task:
$$
\begin{bmatrix}\varphi\\0\end{bmatrix} -  \begin{bmatrix}1\\0\end{bmatrix}
= \begin{bmatrix}\varphi^{-1}\\0 \ \ \ \end{bmatrix}
\,, \qquad
\begin{bmatrix}1\\0\end{bmatrix} - \begin{bmatrix}\varphi^{-1}\\0 \ \ \ \end{bmatrix}
= \begin{bmatrix}\varphi^{-2}\\0 \ \ \  \end{bmatrix}
\,, \qquad \hbox{etc.}
$$


~

~

\newpage
\section{Celestial sphere, modular plane, and Pauli spinors}

In this section we will sent the following observation: 

~\\
{\bf Proposition.}
{\sf    
The celestial sphere image of irreducible Descartes quadruples understood as vectors in the Minkowski space
 coincides with the corresponding lattices
represented in the hyperbolic half-plane.
}

~

Below we explain and prove this observation.

\subsection{Celestial sphere}

\def\m{\phantom{-}}

Descartes' formula may be written in terms of a quadratic form:
\begin{equation}
\label{eq:Q}
\begin{bmatrix}
B_0 \\ B_1 \\ B_2 \\ B_3
\end{bmatrix}^T
\,
\begin{bmatrix}
-1 & \m 1 & \m 1 & \m 1 \\
\m 1 & -1 & \m 1 & \m 1 \\
\m 1 & \m 1 & -1 &\m 1 \\
\m 1 &\m 1 &\m 1 & -1
\end{bmatrix}
\begin{bmatrix}
B_0 \\ B_1 \\ B_2 \\ B_3
\end{bmatrix}
\ \ = \ \ 0
\end{equation}
The matrix $Q$ of this quadratic form defines a pseudo Euclidean norm of signature $(+,-,-,-)$, 
the norm of Minkowski space.%
\footnote{It is however not the Minkowski space of circles as introduced in 
\cite{jk-matrix}, but rather its dual space.}
The above equation implies that the Descartes quadruples 
may be viewed as vectors in the Minkowski space that lie 
on the null cone,
known as the light cone to physicists.%

\begin{figure}[H]
\centering
\includegraphics[scale=.49, angle=-90]{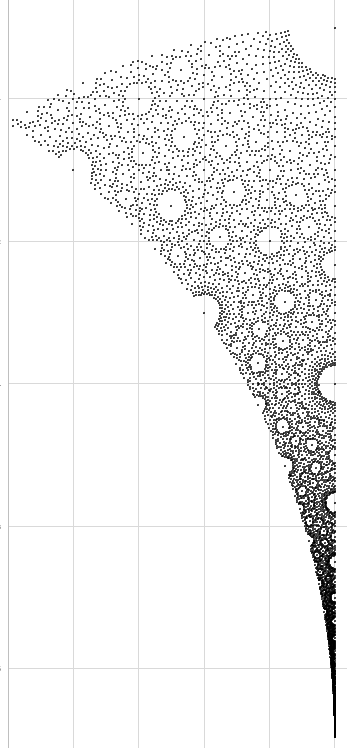}
\caption{Celestial sphere view of Descartes quadruples}
\label{fig:sky1}
\end{figure}

We may think of them as ``photons'' coming from distant stars,
and present them as points in the sky -- the celestial sphere.
A sphere may in turn be stereographically projected on the 2D plane.
The result of this process is presented in Figure \ref{fig:sky1} 
for quadruples with the value of $B$ from 0 to 300
(more than 12000 points).

\begin{figure}[H]
\centering
\includegraphics[scale=.77]{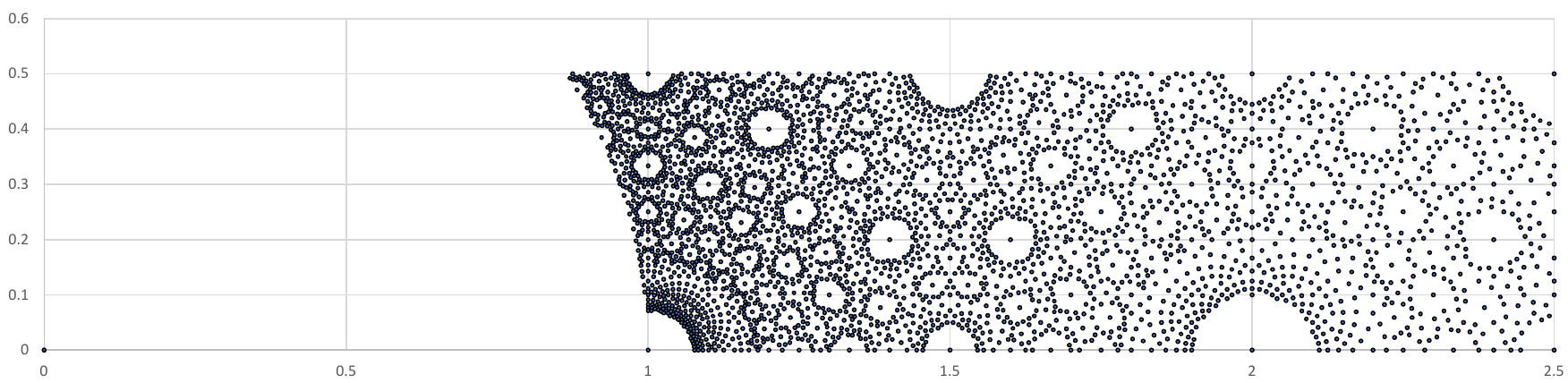}
\caption{Celestial sphere view of Descartes quadruples, projected from the South pole}
\label{fig:sky2}
\end{figure}

Let us trace the steps.
First, diagonalize the quadratic form $Q$ of Equation \eqref{eq:Q}.
Here are two  of many possible diagonalizations of \ref{eq:Q}:
$$
\left\{\begin{array}{rlll}
x&=2(B_2-B_0)   \\
y&=2(B_3-B_1)   \\
z&=2(B_1+B_3) -(B_0+B_2)  \\
t&=2(B_1+B_3) +(B_0+B_2)  
\end{array}\right.
\qquad
\left\{\begin{array}{rlll}
x&=B   \\
y&=\mu   \\
z&=\frac{n-k}{2}  \\
t&=\frac{n+k}{2}  
\end{array}\right.
$$
Descartes formula in the diagonalized coordinates becomes
$$
x^2+y^2+z^2-t^2 \ = \ 0
$$
(as it is easy to check.)
Equation $z=1$ defines a three-dimensional space in $M$, 
and  its intersection with the light cone defines a sphere $S^2$ in it.
It is called the {\bf celestial sphere}.
The Descartes quadruples may be projected along the rays onto this sphere.
$$
x'= x/t,\quad
y' = y/t,\quad
z' = z/t
$$
Obviously, $x'^2 +y'^2 + z'^2 = 1$.
Now, we can perform the stereographic projection from this sphere onto plane via
$$
(X,Y) \ = \ \left(\; \frac{x}{1-z'},\, \frac{y'}{1-z'}\;\right)
\qquad \hbox{or}\qquad
(X,Y) \ = \ \left(\; \frac{x}{1+z'},\, \frac{y'}{1+z'}\;\right)
$$
where the first pair follows the projection from the North pole, while the second from the South pole.
Figure \ref{fig:sky1} shows the latter.
The former would gives an unbounded dust of points shown in Figure \ref{fig:sky2},
extending infinitely to the right. 
\\

\subsection{The modular space of lattices}

There is a standard way to visualize the space of all 2D lattices in a complex plane
as follows (see e.g., \cite{Serre}):
Identify the lattice plane with complex plane.
The basis vectors of a lattice are now two complex numbers $(\omega_1,\,\omega_2)$.
Consider the map:
\begin{equation}
\label{eq:Serre}
(\omega_1,\,\omega_2)\quad \mapsto \quad z = \omega_2/\omega_1
\end{equation}
This map corresponds to scaling and rotating the lattice so that vector $\omega_1$ becomes 1.
Vector $z$ preserves the information of the ratio of the length of the vectors and 
their mutual angle.

\begin{figure}[H]
\centering
\begin{tikzpicture}[scale=2.1]
\clip (-2,0) rectangle (2,2);

\draw [fill=yellow!20, thin] (-1/2,0) rectangle (1/2,2.5);
\draw [fill=white] (0,0) circle (1);
\draw (2.5,0)--(-2.5,0) ;   
\foreach \n in {0, 1, 2, 3}
\draw (\n + 0.5, -1.5)--(\n + 0.5, 2.5) 
          (-\n - 0.5, -1.5)--(-\n - 0.5, 2.5); 
\foreach \n in {-2, -1, 0, 1, 2}
\draw (\n , 0) circle (1/1);
\foreach \n in {-3,-2, -1, 0, 1,2, 3}
\foreach \a/\b   in {  
1/3,  1/5,  1/7,  3/8,  1/9,  1/11,  1/13,  1/15,  4/15,  7/16,  1/17,  1/19,  1/21,  8/21,  1/23,  5/24,  11/24,
1/25,  1/27,  1/29,  1/31,  15/32,  1/33,  1/36,  6/36,  1/37,  1/39,  14/39,  11/40,  19/40,  1/41,  1/43,
1/45,  19/45,  1/47,  7/48,  23/48,  1/49,  1/51,  16/51,  1/53,  1/55,  21/55,  13/56,  27/56,  1/57,  20/57,
1/59, 1/61, 1/63, 8/63, 31/64, 1/65, 14/65, 1/67, 1/69, 22/69 
}
\draw (\n + \a/\b, 0) circle (1/\b)
           (\n - \a/\b, 0) circle (1/\b)
;
\end{tikzpicture}
\caption{The Dedekind tessellation}
\label{fig:wiki}
\end{figure}
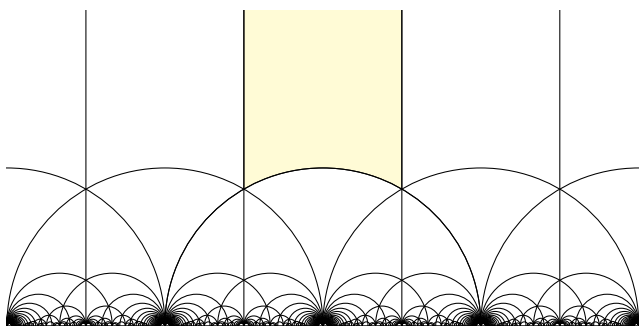

If one applies this to the principal basis, 
its complex representation $z$ lies in the fundamental domain
of the modular tessellation of the Poincar{\'e} half-plane (the shaded region in Figure \ref{fig:wiki}),
more precisely, in the region 
defined by $|z|\geq 1$,\; $\hbox{\rm Im}\,z >0$, \; $0\leq\hbox{\rm Re}(z) < 1/2$.
The other bases fall in the other tiles of the Dedekind tessellation
and are images of the principal point $z$ via M\"obius action 
$$
z \ \mapsto \ \begin{bmatrix} a & b \\ c & d\end{bmatrix} \cdot z = \frac{az+b}{cz+d}
$$
where the matrix is a group element of the unimodular group $\SL(2,\mathbb Z)$.

Let us apply this to the lattices representing the Descartes configurations
with the principal spinors $\mathbf a$ and $\mathbf b$: 
$$
\begin{array}{rcl}
z &=& \dfrac{a}{b} \ = \ \dfrac{|\mathbf b|}{|\mathbf a|} \;{\rm e}^{i\theta}\\[12pt]
  &=& \dfrac{|\mathbf b|}{|\mathbf a|} \cos\theta +   i  \,  \dfrac{|\mathbf b|}{|\mathbf a|} \sin\theta  \\[12pt]
  &=& \dfrac{|\mathbf a| \/ |\mathbf b|\/ \cos\varphi}{\|\mathbf a\|^2}  
               + \dfrac{|\mathbf a| \/ |\mathbf b|\/ \sin\varphi}{\|\mathbf a\|^2}   \, i  \\[12pt]
  &=& \dfrac{|\mathbf a \bigcdot \mathbf b|}{\|\mathbf a\|^2}  
               + \dfrac{|\mathbf a\times \mathbf b|}{\|\mathbf a\|^2}   \, i 
\end{array}
$$
Thus
\begin{equation}
\label{eq:z-k}
z = \dfrac{\mu}{k}  + \dfrac{B}{k}   \, i  
 \end{equation}
With the help of software, we can produce an image that turns out quite similar to the one in Figure
\ref{fig:sky2}.

\subsection{Proof of the Proposition 5.1}

The puzzle of similarity of the patterns in the two images, 
the one obtained via celestial sphere and one in the hyperbolic plane,
may easily be explained. 
One of the diagonalizations of the Descartes formula is:
$$
B^2 + \mu^2 + \left(\frac{n-k}{2}\right)^2 - \left(\frac{n+k}{2}\right)^2 \ = \ 0
$$
Associate  $x=B$, $y = \mu$, $z= \frac{n-k}{2}$, $t= \frac{n+k}{2}$.
Substituting these entities to the celestial sphere coordinates gives:
$$
x' = \frac{x}{t} = \frac{2B}{n+k} , \qquad
y' = \frac{y}{t} = \frac {2\mu}{n+k}, \qquad
z' = \frac{n-k}{n+k}
$$
and consequently
\begin{equation}
\label{eq:s-k}
X \ = \ \frac{x'}{1- z'} \ = \  \frac{B}{k} , \qquad
Y  \ = \ \frac{y'}{1-z'} \ = \ \frac{\mu}{k} 
\end{equation}
which coincides with the formula for the projection from the North projection of the celestial sphere \eqref{eq:z-k}.
Hence the obtained image must coincide with the belt-like right half of the fundamental region 
of Figure \ref{fig:wiki}.

The stereographic projection of the celestial sphere from the south pole S requires 
a sum in the denominator of \eqref{eq:s-k} and resolves as follows:
\begin{equation}
\label{eq:s-n}
X \ = \ \frac{x'}{1+ z'} \ = \  \frac{B}{n} , \qquad
Y  \ = \ \frac{y'}{1+z'} \ = \ \frac{\mu}{n} 
\end{equation}
($n$ instead of $k$ in the denominator.)
The resulting image is illustrated in Figure \ref{fig:sky1}.
Now, for the complex plane representation of lattices,
one obtains the matching image by mapping the fundamental domain 
to the triangle-like region underneath
by the unimodular transformation:
$$
z\mapsto z' \ = \  \begin{bmatrix}0&-1\\ 1&0\end{bmatrix}\cdot z \ = \ \frac{-1}{\dfrac{\mu}{k}  + \dfrac{B}{k}   \, i} 
\ = \ 
\frac{-k}{\mu + B  \, i} 
\ = \ 
\frac{-k(\mu - B  \, i)} {\mu^2+B^2}
\ = \ 
\frac{-\mu + B  \, i} {n}
$$
(after using $B^2+\mu^2=kn$).    
Reflection will map the point to the right side of the triangular region, with coordinates
$$
z' = \frac{\mu}{n} + \frac{B}  {n}\, i
$$
which coincides with the formulas for  projection of the celestial sphere projected onto the plane from the south,
\eqref{eq:s-n}.

\begin{figure}[H]
\centering
\includegraphics[scale=.32]{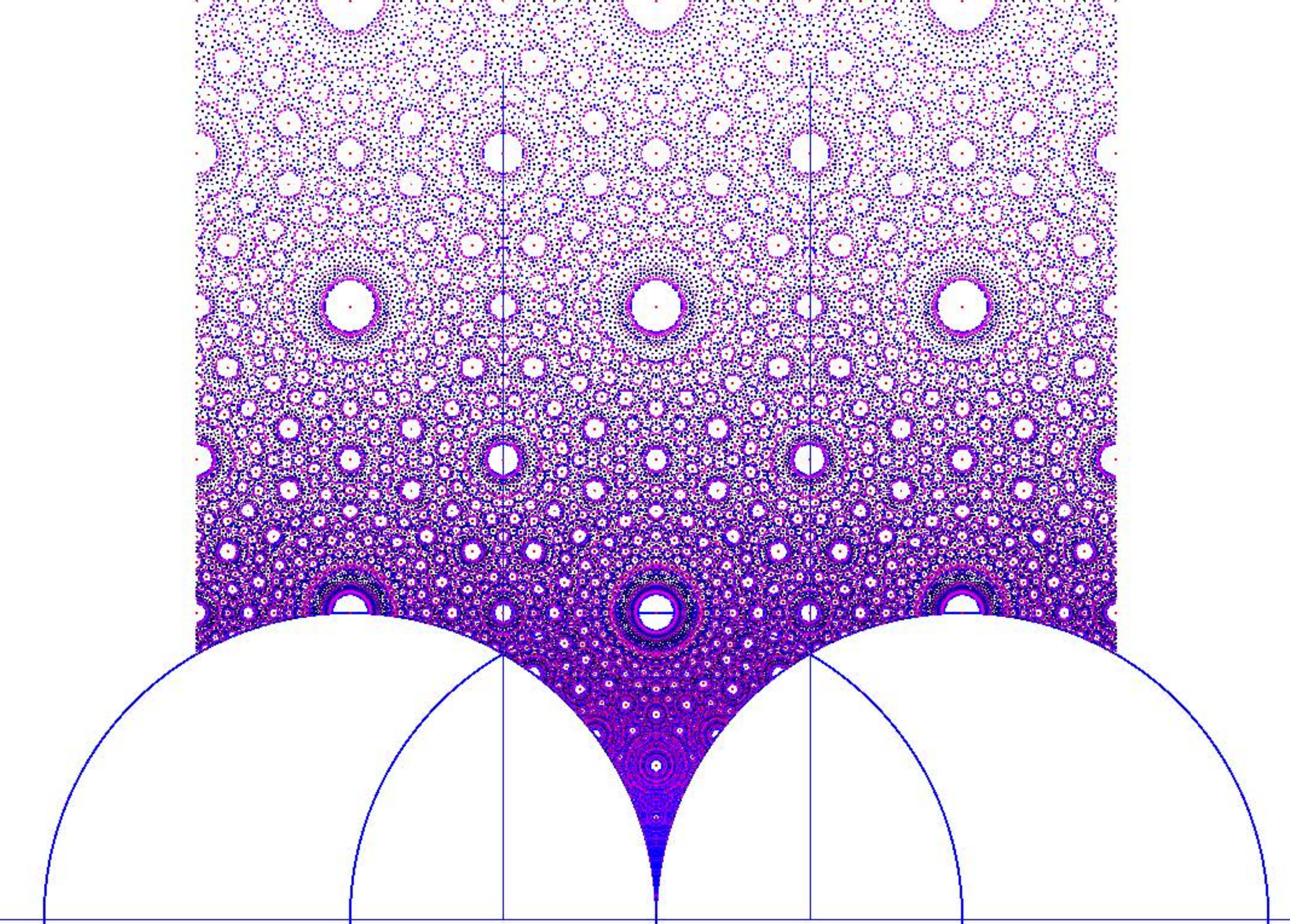}
\caption{Descartes dust in the fundamental domain and its images in a few adjacent regions}
\label{fig:Ddust}
\end{figure}

\newpage

\subsection{Unification: Pauli spinors from tangency spinors}

The tangency spinors are spinors for the 3-dimensional space with Minkowski
quadratic form of signature (1,2).
It has been a somewhat metaphorical use of the term native to mathematical physics
of (1+3)-dimensional spacetime
for objects in 2-dimensional geometry.
But now, by some magical turn, the two worlds merge
when we represent a pair of tangency spinors
as a single Pauli spinor:

~\\
Consider the following map
$
\mathbb R^2 \times \mathbb R^2 \ \to \ \mathbb C^2
$.
First, view the tangency spinors back as the complex numbers:  
$\mathbf a = \begin{bmatrix}a_1\\ a_2\end{bmatrix}\ \mapsto \ a = a_1+ia_2$
(same for spinor  $\mathbf b$).
Next, define a ``metaspinor''
$$
\mathbf u = \begin{bmatrix} a \\ b\end{bmatrix} \ \in \mathbb C^2
$$ 
(which is the usual Pauli spinor,  known in the context of (1+3)-Minkowski space).
Now, consider a Hermitian matrix constructed from this spinor:
$$
H \ = \ \mathbf u \otimes \mathbf u^*
$$
where the star denotes the Hermitian conjugation: transposition and complex conjugation.
The tensor product may be realized as the Kronecker product.

~

\noindent
{\bf Proposition:}
{\sf
The connection between the tangency spinors and the curvatures,
 and the parametrizing formula (and consequently the Descartes formula)
  can be expressed in a single statement:
$$
\begin{bmatrix}
   k & \mu -iB\\
   \mu +iB & n \end{bmatrix}
   \ = \ 
\begin{bmatrix} \mathbf a  \\ \mathbf b \end{bmatrix}
\otimes
\begin{bmatrix} \mathbf a \\ \mathbf b \end{bmatrix}^*
\qquad  \xrightarrow{\quad\rm det \quad} \ 0 \  
$$
}

\noindent
{\bf Proof:}
We get 
$$
\begin{array}{rcl}
H&=& \mathbf u \otimes \mathbf u^* 
\quad = \quad
\begin{bmatrix} a \\ b \end{bmatrix} 
    \begin{bmatrix} \overline{a} & \overline{b} \end{bmatrix} 
\quad = \quad 
\begin{bmatrix} |a|^2 & a\overline b\\
                           b\overline a & |b|^2 \end{bmatrix}                \\[15pt]
&=&  \begin{bmatrix} |a|^2 & a\!\bigcdot\! b - (a\!\times\! b)\, i \\
                           a\!\bigcdot\! b   +  (a\!\times\! b)\, i & |b|^2 \end{bmatrix} \\[15pt]
&=&    \begin{bmatrix} k & \mu-iB\\
                           \mu + iB & n\end{bmatrix} 
                           \end{array}
$$
Since the matrix $\mathbf u \otimes \mathbf u^* $ is  of rank 1, its determinant vanishes automatically.
This immediately implies $B^2+\mu^2 =kn$.
\QED

~\\
{\bf Remark 1:} The above parallels the correspondence well known in theoretical physics,
where vectors $\mathbf v$ in Minkowski space-time are mapped to the space of Hermitian matrices:
 $$
\mathbf v = [x,y,z,t] \quad \mapsto \quad M=\begin{bmatrix}t-z & y-ix\\
                                                    y+ix & t+z   \end{bmatrix}
$$
The Lorentz norm is represented now by a determinant:
$$
\|\mathbf v\| \ = \ t^2 - x^2 - y^2 - z^2 \ = \ \det M
$$
In particular, any vector (Pauli spinor) $\mathbf u = [u_1\ u_2]^T\in \mathbb C^2$ defines 
a null-vector (``photon'')
by setting $M=\mathbf u\otimes \mathbf u^*$,
which may be understood as the Kronecker product. Since it is of rank 1, its determinant is vanishing,
and so is the norm of the Minkowski vector it defines.
Once we view Descartes quadruples as null vectors in the Minkowski space,
we may evoke these standard means from physics,
including the Pauli spinors, the elements of $\mathbb C$. 
This is where one may use Pauli matrices:
$$
\sigma_0 \ = \  \begin{bmatrix} 1 & 0 \\ 0 & 1 \end{bmatrix}
\qquad
\sigma_1 \ = \  \begin{bmatrix} 0 & 1 \\ 1 & 0 \end{bmatrix}
\qquad
\sigma_2 \ = \  \begin{bmatrix} 0 & -i \\ i & 0 \end{bmatrix}
\qquad
\sigma_3 \ = \  \begin{bmatrix} 1 & 0 \\ 0 & -1 \end{bmatrix}
$$
which allow one to express a vector as a Hermitian matrix via the map:
$\mathbf v \ \mapsto \ t\sigma_0 + x\sigma_1 + y\sigma_2 + z \sigma_3$

~\\
{\bf Remark 2:}
Some additional formulas for matrix $H$ follow easily:
$$
\begin{array}{rcl}
H 
%
&=& \begin{bmatrix}  B_0+B_1 & \sqrt{B_1B_2\!+\!B_2B_3\!+\!B_3B_1} + i B_0 \; \\
                           \; \sqrt{B_1B_2\!+\!B_2B_3\!+\!B_3B_1} - i B_0 & B_0+B_2 \phantom{\Big|}
                           \end{bmatrix} \\[23pt]
&=& \begin{bmatrix} B_0+B_1 & \frac{1}{2}(B_0\! +\!B_1\!+\!B_2\!-\!B_3) + i B_0 \; \\
                           \;\frac{1}{2}(B_0\! +\!B_1\!+\!B_2\!-\!B_3) - i B_0 & B_0+B_2 \end{bmatrix} 
\end{array}                           
$$

\subsection{Hopf fibration and the projective spinor space}

The stardust may also be reproduced in terms of the Pauli spinors.
Indeed,
consider the following chain of maps
$$
\mathbb C^2  
\ \  \xrightarrow{~\text{norm}~} \ \ 
S^3  
\ \ \xrightarrow{~\text{Hopf}~} \ \
S^2 
\ \  \xrightarrow{~\text{stereo}~} \ 
\mathbb R^2 
$$
The first map in the chain is a normalization of vectors by the standard Hermitian norm in $\mathbb C^2$.
The result is a three-dimensional sphere in $\mathbf C^2$.
The second map is the Hopf fibration of the 3D sphere over 2D sphere $S^2$,
The composition of these maps is simply the projection map from $\mathbb C^2$
to the projective space $\hbox{P}\mathbb C^2 \cong S^2$. 
In such a context, $S^2$ is known in quantum physics as the Bloch sphere.
The last map is the usual stereographic projection of a 2D sphere to the plane. 
As it is easy to see,
the composition of these maps abbreviates to
$$
\begin{bmatrix}a\\b\end{bmatrix} 
\ \mapsto \ 
\begin{bmatrix}1\\b/a\end{bmatrix}
\ \mapsto \ z = b/a
$$
which is exactly of the same form as the maps considered in the previous subsections,
cf., \eqref{eq:Serre}.
(Clearly, the zero vector needs  to be excluded from these maps.)


\begin{figure}[H]
\hspace{-.2in}%
\includegraphics[scale=.74]{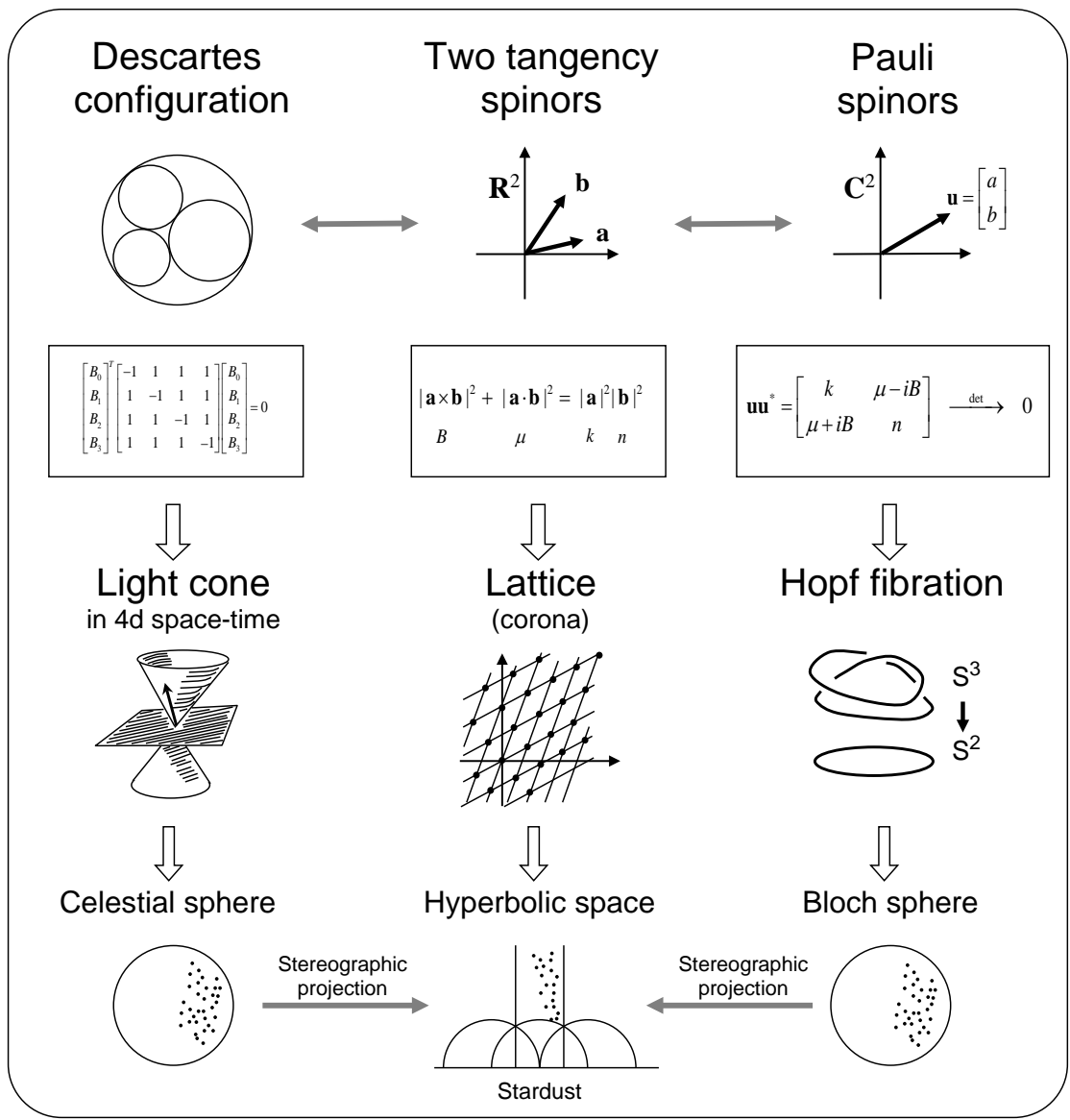}
\caption{Main ingredients of the relations.  
Various forms of the Descartes formula are framed.}
\label{fig:diagram}
\end{figure}
%
%
%
%
%
%
%
%
%
%

\newpage

\section{Summary}


\noindent
{\bf 1.}  
The relations between the tangency spinors and Descartes configurations
discussed in the present paper are summarized in the following diagram.  
The top row  concerns spinors, the bottom the disks. 
Restricting the objects to integral curvatures 
leads to the parametrization of the irreducible packings. 
\def\smalll{\scriptsize}
%
$$
\begin{tikzpicture}[baseline=-0.8ex]
    \matrix (m) [ matrix of math nodes,
                         row sep=2.8em,
                         column sep=4em,
                         text height=3.8ex, text depth=3ex] 
   {
   \quad { \hbox{\sf arbitrary} \atop  {\hbox{\sf two 2D} \atop \hbox{vectors}} } \quad  
      & \quad  { {\hbox{\sf Lattice} \atop \hbox{\sf (mosaic)}} } \quad  
      & \quad {\hbox{\sf Principal}\atop\hbox{\sf basis}} \quad \\
   \quad{\hbox{\sf everted} \atop \hbox{tricycle}}  \quad  
      & \quad {\hbox{\sf Apollonian packing} \atop \hbox{\sf (corona)}} \quad  
      & \quad {\hbox{\sf maximal}\atop\hbox{\sf tricycle}} \quad \\
     };
    \path[<->]
        (m-1-1)   edge node[above] {\sf\smalll } (m-1-2)
        (m-1-2)   edge  (m-1-3)
        (m-2-1)   edge   (m-2-2)
        (m-2-2) edge node[above]  {\sf\smalll  } node[below]  {$\hbox{\sf\smalll generate}$}(m-2-3)
        (m-1-1) edge node[right] {\sf\smalll  } (m-2-1)
        (m-1-2) edge node[right] {\sf\smalll  } (m-2-2)
        (m-1-3) edge node[right] {\sf\smalll  } (m-2-3)
        ;
\end{tikzpicture}   
$$
%
Two arbitrary 2D integral vectors $\mathbf v$ and $\mathbf w$
define an integral tricycle, which, when extended, results in 
an integral Apollonian disk packing $\mathcal A(\mathbf v, \mathbf w)$.
But they also define a lattice $L(\mathbf v, \mathbf w)$ in the spinor space via linear combinations.
One may find the principal basis $(\mathbf a, \mathbf b)$ of this lattice;
clearly $L(\mathbf a, \mathbf b)=L(\mathbf v, \mathbf w)$. 
There are two benefits of this finding:
\vspace{-.05in}
\begin{itemize}
\setlength\itemsep{0.14em}%
 \item[1.] 
Vectors $(\mathbf a, \mathbf b)$ relate directly to the terms 
of the parametrization theorem for Apollonian disk packings 
and explain their meaning.
\item[2.]
The principal spinors point to the maximal tricycle (and maximal Descartes configuration)in the Apollonian packing.
\end{itemize}
\vspace{-.05in}
\noindent
Other circles in the boundary corona in $\mathcal A$ 
are obtainable via linear combinations of the basis vectors.  
\\

\noindent
{\bf 2.}  
Somewhat unexpectedly, 
there is a deeper level on which the concepts are related.
The results may be reorganized 
with the help 
of a number of ideas known from theoretical physics,
like Pauli spinors of quantum mechanics 
or celestial sphere from the theory of relativity.
Also, the family of the integral Apollonian packings finds
a geometric image as points
in celestial sphere, Bloch sphere, or hyperbolic space.  
This is summarized in Figure~\ref{fig:diagram}.

~

\noindent
{\bf 3.}  
On the lighter side, as a final remark, let us notice that   
the two main spinor equations lead to two different geometric extensions:
$$
\begin{array}{rlccc} 
\curlrm \mathbf u = 0 & \qquad \to \quad \hbox{tessellations}\\[2pt]
\divrm \mathbf u = 0  & \qquad \to \quad \hbox{lattices}
\end{array}
$$
The first connection is explained in  \cite{jk-tessellation}. Stripping its geometric content,  we get 
a method of obtaining integral Descartes configurations (and consequently integral Apollonian disk packings) 
from four arbitrary integers as follows.

\newpage


\noindent
{$\bullet$} 
Let $\mathbf a$ and $\mathbf b$ be arbitrary integral vectors in $\mathbb Z^2$.  Then the following 
values satisfy the Descartes equation \eqref{eq:Descartes} for three tangent inner disks
and the fourth one of the two conjugated disks $D_1$ or $D_2$:
\begin{equation}
\label{eq:reduced}
\begin{array}{rl}
A  &= \ \  |\mathbf b|^2 +\mathbf a\bigcdot \mathbf b\\
B  &=\  \   |\mathbf a|^2 +\mathbf a\bigcdot \mathbf b \\
C  &=\ \ -\mathbf a\bigcdot \mathbf b\\
D_{1,2}  &=  \ \ |\mathbf a|^2+|\mathbf b|^2  + \mathbf a\bigcdot \mathbf b 
\; \pm \; 2\, 
\mathbf a\times \mathbf b
\end{array}
\end{equation}
%
{$\bullet$} Compare it with the content of the present paper:
\begin{equation}
\label{eq:reduced2}
\begin{array}{rl}
A  &= \ \  |\mathbf a|^2 +|\mathbf a\times \mathbf b|\\
B  &=\  \   |\mathbf b|^2 +|\mathbf a\times \mathbf b| \\
C  &=\ \ -|\mathbf a\times \mathbf  b|\\
D_{1,2}  &=  \ \ |\mathbf a|^2+|\mathbf b|^2  \pm 2 \mathbf a\bigcdot \mathbf b 
+ |\mathbf a\times \mathbf b|
\end{array}
\end{equation}

On the other hand,
the spinor - Descartes quadruple can also be given a simple 
geometric interpretation in terms of areas of tiles,
as shown in Figure \ref{fig:tiles}.
The parallelogram is that of the fundamental cell.
The curvatures of the greatest disk in the implied Descartes configuration
are the sums of the dark (red) region $B$ and one of the two lighter squares (yellow).
The next two (shown aside) are based on the two diagonals of the parallelogram.

\def\nmaly{.47}
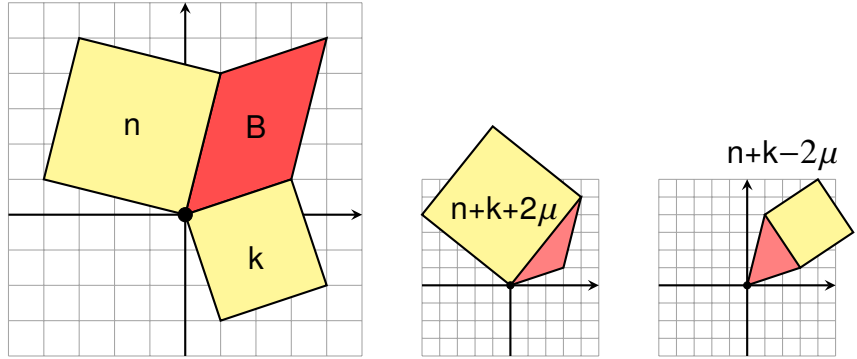
\begin{figure}[H]
\centering
$$
\begin{tikzpicture}[scale=\nmaly, rotate=0, shift={(0,0cm)}]  
\draw[thin, color=gray!70] (-5,-4) grid (5,6);
\draw[-stealth,  thick] (-5,0) -- (5,0);   
\draw[-stealth,  thick] (0, -4) -- (0, 6);   
\draw [thick, fill=red!70,fill opacity=.7]  (0,0) -- (3,1) --  (3+1, 1+4)  -- (1,4) -- cycle;

\draw [thick, fill=yellow!50, fill opacity=.7]  (0,0)  -- (1,4) -- ( -3,5) --(-4,1)-- cycle;
\draw [thick, fill=yellow!50, fill opacity=.7]  (0,0)  -- (3,1) -- ( 4,-2) --(1,-3)-- cycle;

\draw node at  (2,2.5) [scale=1.2]{\sf B};
\draw node at  (2,-1.2) [scale=1.2]{\sf k};
\draw node at  (-1.5, 2.5) [scale=1.2]{\sf n};

\draw [fill=black] (0,0) circle  (.2);
\end{tikzpicture}%
\qquad 
\begin{tikzpicture}[scale=.5*\nmaly, rotate=0, shift={(0,0cm)}]  
\draw[thin, color=gray!70] (-5,-4) grid (5,6);
\draw[-stealth,  thick] (-5,0) -- (5,0);   
\draw[-stealth,  thick] (0, -4) -- (0, 6);   
\draw [thick, fill=red!50, fill opacity=.7]  (0,0) -- (3,1) --  (3+1, 1+4)  -- (1,4) -- cycle;

\draw [thick, fill=yellow!50, fill opacity=.7]  (0,0)  -- (4,5) -- ( -1,9) --(-5,4)-- cycle;

\draw node at  (-.35,4.2) [scale=1.2]{\sf n+k+2$\mu$};

\draw [fill=black] (0,0) circle  (.2);
\end{tikzpicture}%
\qquad
\begin{tikzpicture}[scale=.5*\nmaly, rotate=0, shift={(0,0cm)}]  
\draw[thin, color=gray!70] (-5,-4) grid (5,6);
\draw[-stealth,  thick] (-5,0) -- (5,0);   
\draw[-stealth,  thick] (0, -4) -- (0, 6);   
\draw [thick, fill=red!50, fill opacity=.7]  (0,0) -- (3,1) --  (3+1, 1+4)  -- (1,4) -- cycle;

\draw [thick, fill=yellow!50, fill opacity=.7]  (1,4)  -- (3,1) -- ( 6,3) --(4,6)-- cycle;

\draw node at  (2,7.4) [scale=1.2]{\sf n+k$-2\mu$};

\draw [fill=black] (0,0) circle  (.2);
\end{tikzpicture}%
$$

\caption{Curvatures of Descartes configuration via simple geometry}
\label{fig:tiles}
\end{figure}

\noindent
{\bf Note on software:}
Most of the figures were made with Tikz \cite{tikz}.

\section*{Appendix A}

\begin{figure}[H]
\centering
\includegraphics[scale=.85]{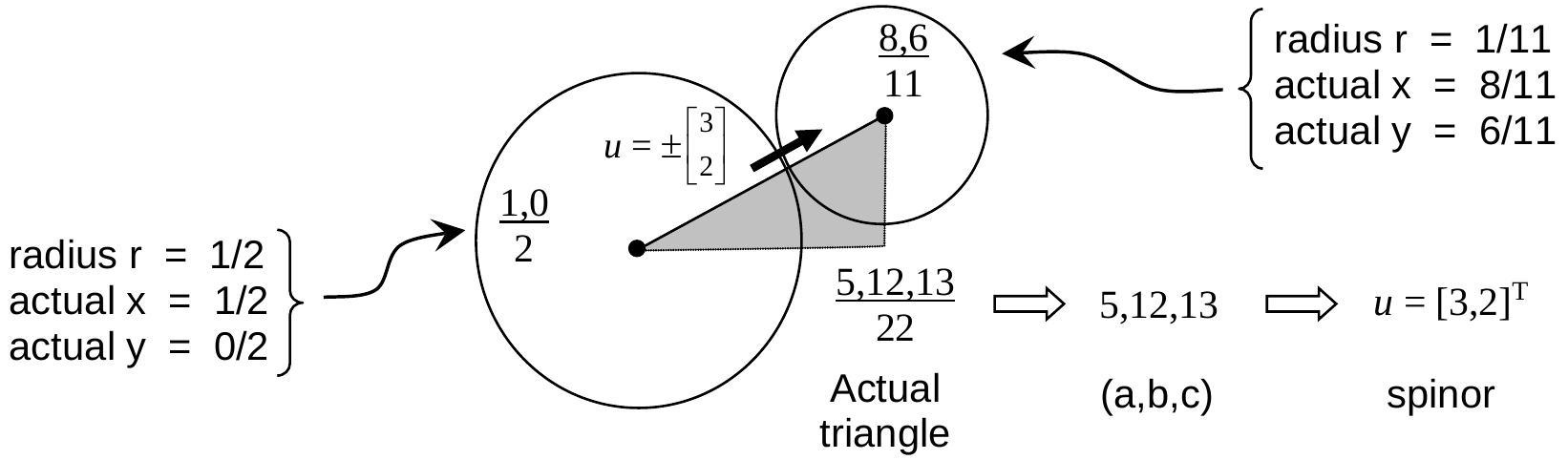}
\caption{\small From two tangent circles to a spinor (not to scale)}
\label{fig:process}
\end{figure}

The geometric interpretation and motivation follows.
Every disk in the Cartesian plane may be given a {\bf symbol},
a fraction-like label that encodes the size and position of the disk:  
the curvature  is indicated in the denominator 
while the positions of the centers may be read off by interpreting the symbol as a pair of fractions
\cite{jk-matrix}.
$$
\hbox{\sf symbol:\  } 
       \frac{\dot x,\;\dot y}{\beta}  
                         \qquad \Longrightarrow \qquad 
       \begin{cases}     \hbox{\sf radius:\  } &  r = \frac{1}{\beta} \\
                                  \hbox{\sf center:\ }  & (x,\,y)  =   \left(\frac{\dot x}{\beta},\, \frac{\dot y}{\beta}\right) 
       \end{cases}
$$
The numerator, called the {\bf reduced coordinates} of the a disk's center, 
is denoted by dotted letters $(\dot x, \, \dot y)= (x/r,\, y/r)$.
Unbounded disks extending outside a circle are given negative radius and curvature.
Two tangent disks in a plane  define a triangle
with sides as follows \cite{jk-Clifford,jk-spinor}:
\begin{equation}
\label{eq:product}
     \frac{\dot x_1, \;  \dot y_1}{ \beta_1}
\ \Join  \
     \frac{\dot x_2,\; \dot y_2}{ \beta_2}
\qquad \mapsto \qquad
       \left[\begin{matrix}
                     a \\
                     b \\
                    c
     \end{matrix}\right]
\equiv 
  \left[\begin{matrix}
                     \beta_1 \dot x_2 - \beta_2 \dot x_1 \\
                     \beta_1 \dot y_2 - \beta_2 \dot y_1 \\
                    \beta_1+\beta_2
     \end{matrix}\right]
\end{equation}
where $a^2+b^2=c^2$ 
(see Figure \ref{fig:process}).
The actual size of the triangle in the plane is scaled down by the factor of $\beta_1\beta_2$
(gray triangles in Figure~\ref{fig:process}).
The symbols in some disk packing are integral, then so are the triples $(a,b,c)$.  
Recall that Pythagorean triangles admit Euclidean parameters 
that determine them via the following prescription:
$$
\mathbf u=\vvec mn \quad \to \quad (a,b,c) = (m^2-n^2,\ 2mn, \ m^2+n^2)
$$
(see, e.g., \cite{Sie, T-T}).
As explained in  \cite{jk-Clifford}, Euclidean parameters can be viewed as a {\it spinor},
a vector of $\mathbb T \cong \mathbf u\in \mathbb R^2$.
Equivalently, viewing the spinor as a complex number $u\in\mathbb C\cong\mathbb R^2$
the above relations is defined by squaring:
$$
            u=
            m+ni \quad \to \quad u^2  = a+bi   = (m^2-n^2) + 2mn\, i 
$$
with $c=|u^2|=m^2+n^2$.  
We extend  this map to arbitrary oriented triangles, not necessarily integer.
\\

The emergence of the {\bf tangency spinor} for a pair of tangent disks  
is summarized in Figure~\ref{fig:process}.

%

The curvatures are the numbers in the denominators of the ``symbols'' inscribed into the disks.
The numerator numbers in the fractions code the positions of the disks in the euclidean plane
and  are of no concern for us here; they are explained in Appendix A.

\subsection*{Drawing the Apollonian packing from the defining spinors}

It is easy to define uniquely the spinor given full data on the circles.
But the reverse process can be done up to translations.

~

\noindent
{\bf Proposition A.1:}
{\sf Let us assume that the main circle, the boundary of the outer disk of negative curvature $B_0=-B$
is located at the origin:
$$
\hbox{center at} \ (0,0) \qquad \hbox{radius = } 1/B
$$
Let $A$ be a disk in the major corona, the corona of $B_0$,
and let $\mathbf a$ be a spinor
$$
\mathbf a=\begin{bmatrix}m \\ n\end{bmatrix}  =  \spin(A,B_0)
$$
(in the graphical representation it is drawn as if leaving the packing),
Then the position and the radius of the disk $A$ are: 
$$
\hbox{center at} \ \left(\frac{m^2-n^2}{B(B+m^2+n^2)},\frac{2mn}{B(B+m^2+n^2)}\right) 
\qquad 
\hbox{radius = } 1/A
$$
In terms of the disk symbol:
$$
\frac{(m^2-n^2)/B,\;  2mn/B}{m^2+n^2+B} 
$$
}

\noindent
{\bf Proof:} This is a simple consequence of the definition of the spinor $\spin (K,A)$
for two disks of radii $r_K=1/K$ and $r_A=1/A$, respectively.
If viewed as a complex number, it defines
the vector distance between the centers of two circles $K$ and $A$
as $\overrightarrow{AK} = u^2r_A r_K$.
Hence, if the position of disk $K$ is $p\in\mathbb C$, the position of circle $A$ is 
$p+u^2/(KA)$. 
A difference instead of sum must  be used for the opposite spin $\spin(AK)$.
Apply it to the situation described in the proposition to get the result.
\QED

~

Another question is how to construct spinors given only curvatures
of the disks.

~

\noindent
{\bf Proposition:}
For any choice of three curvatures of three tangent disks,
$(B_0, B_1, B_2)$,  with $B_0<0$, there exists a spinor description, namely
$$
\mathbf a = \begin{bmatrix} \dfrac{B_0}{\sqrt{B_0+B_1}}  \\ \vspace{-.07in} \\ 
                                             \dfrac{\strut\sqrt{B_0B_1+B_1B_2+B_2B_0}}{\sqrt{B_0+B_1}}  \end{bmatrix} 
                                             \qquad
\mathbf b = \begin{bmatrix} \phantom{\dfrac{A}{\sqrt{B}}} 0 \phantom{\dfrac{A}{\sqrt{B}}}  \\ \vspace{-.07in} \\ 
                                              \dfrac{B_0+B_1}{\sqrt{B_0+B_1}}\end{bmatrix}  
$$
The numerator in the second entry of $\mathbf a$ may be replaced by $B_0+B+1+B_2-B_3)/2$,
where $B_3$ is a solution to the Descartes problem for the given three disks. 
Its value coincides with the curvature of mid-circle $(B_0,B_1,B_2)$.

\section*{Appendix B: the symmetry group aspect}

A few words on the group symmetries that pertain to the Descartes configuration 
and the corresponding spinor structure are in order.
The mutual relations are displayed in the diagrams that follow.
The Lorentz group $\SO_{(1,3)}(\mathbb R)$ acting on the Descartes configurations will 
permute them. 
The discrete subgroup $\SO_{(1,3)}(\mathbb Z)$ carries integral Descartes configurations to likewise ones.   
The Pauli spinors transform by the corresponding 
elements of the group $\SL(2,\mathbb C)$.

The first diagram has a similar content, except it concerns
a single tangency spinor originating from a Pythagorean triangle.
The dashed line in the first diagram indicates the Euclidean parametrization.
In the second, the Pauli spinor describing a Descartes configuration.

~\\
{\bf Pythagorean triples and Euclid's parametrization}
\\
\\
\begin{tikzpicture}[baseline=-0.8ex]
    \matrix (m) [ matrix of math nodes,
                         row sep=2em,
                         column sep=3.5em,
                         text height=3.8ex, text depth=3ex] 
      {
            \dgaction{\ \ \SO(1,2)\ }{\mathbb R^{1,2}}  \quad   
                & \quad \dggaction{\ \ \SL(2,\mathbb R) \  }{\ \ \hbox{Symm}_\circ(2, \mathbb R)\ \ }  \quad 
                & \quad \dgaction{\ \ \SL(2,\mathbb R) \  }{  \mathbb R^2 } \quad  \\
            \dgaction{\ \ \SO(1,2)\  }{\hbox{\sf\smalll  Light Cone}} \quad
                & \quad \dggaction{\ \ \SL(2,\mathbb R) \  }{\ \ \hbox{Symm}_\circ^\bullet(2, \mathbb R) \ \ }       \quad
                & \quad \dggaction{\ \ \SL(2,\mathbb R) \  }{\ {\mathbb R}^2\otimes{\mathbb R}^2 \ } \quad       \\
            \dgaction{\ \ \SO(1,2)\  }{S^1} \quad
                & {}
                &{}&[-1.75cm] \quad \dgaction{\ \ \SL(2,\mathbb R) \  }
                                     { \ {\Prm^2}\mathbb R  \equiv \dot{\mathbb R}\equiv \mathbb R\!\cup\!\{\infty\}\    } \quad       \\
     };

    \path[-stealth, very thick]
        (m-1-2) edge node[above ] { \sf\smalll  2:1}  (m-1-1)
       (m-2-3) edge node[above] {$\hbox{ \sf\smalll  1:1 }$}    
                             node[below] {$\hbox{ \sf\smalll  2:1 }$}   (m-2-2)
       (m-2-2) edge [dotted] node[above] {$\hbox{ \sf\smalll  2:1 }$}   
                                           node [below] {\sf\smalll 1:1}  (m-2-1)
       (m-3-1) edge [dotted,-, very thick] 
                                           node [above] {\sf\smalll stereographic projection} 
                                           node [below] {\sf\smalll (angle doubling)} 
                                           (m-3-4)

        (m-1-1) edge  node[right]  {$\imath$} 
                             node[color=red, above=50pt] 
                                                                 {\smalll\sf Triangles~~~ }                                  (m-2-1)
        (m-2-1) edge node[right] {$\pi $}                                     (m-3-1)
        (m-2-2) edge node[right] {$\imath$}                                     (m-1-2)
           (m-1-3) edge [double,thick]  node[right] {$\hbox{id}\otimes *$} 
                                                     node[color=red, above=50pt] 
                                                                    {\smalll\sf Euclid's parameters }            (m-2-3)
           (m-1-3.east) edge [dotted,bend left =35pt, out=60, in=150]  node[right] {$\pi$}                     (m-3-4)
           (m-1-3) edge [dashed, thin] node[above, rotate=25] 
                                    {$\hbox{\sf\smalll Euclid~} \atop \hbox{\sf\smalll parameters}$  }                 (m-2-2)
;
\end{tikzpicture}   

~\\
{\bf Physics, relativity, space-time}
\\\\
\begin{tikzpicture}[baseline=-0.8ex]
    \matrix (m) [ matrix of math nodes,
                         row sep=2em,
                         column sep=3.5em,
                         text height=3.8ex, text depth=3ex] 
      {
            \dgaction{\ \ \SO(1,3)\ }{\mathbb R^{1,3}}  \quad   
                & \quad \dggaction{\ \ \SL(2,\mathbb C) \  }{ \ \ \hbox{Herm}(2, \mathbb C)\ \ }  \quad 
                & \quad \dgaction{\ \ \SL(2,\mathbb C) \  }{  \mathbb C^2 } \quad  \\
            \dgaction{\ \ \SO(1,3)\  }{\ \hbox{\sf\smalll Light Cone}\ } \quad
                & \quad \dggaction{\ \ \SL(2,\mathbb C) \  }{\ \ \hbox{Herm}^\bullet(2,\mathbb C) \ \ }       \quad
                & \quad \dggaction{\ \ \SL(2,\mathbb C) \  }{\ \ {\mathbb C}^2\otimes{\mathbb C}^2 \ \ } \quad       \\
            \dgaction{\ \ \SO(1,3)\  }{ S^2 \ \ } \quad   
                & {}
                &{}&[-1.95cm] \quad \dgaction{\ \SL(2,\mathbb C) }
                            {\ \ {\Prm^2}\mathbb C  \equiv 
                           {  \dot{\mathbb C}}\equiv \mathbb R^2\!\cup\!\{\infty\}\phantom{\big|} } \quad       \\
     };

    \path[-stealth, thick]
        (m-1-2) edge node[above ] { \sf\smalll  2:1} 
                              node[below ] { \sf\smalll  1:1}    
                              node[above=21pt, color=red ] {\sf\smalll Dirac's  belt trick}    (m-1-1)
       (m-2-3) edge node[above] {$\hbox{ \sf\smalll  induce }$}     (m-2-2)
       (m-2-2) edge [dotted] node[above] {$\hbox{ \sf\smalll  2:1 }$}   
                                           node [below] {\sf\smalll }  (m-2-1)
       (m-3-1) edge [dotted,-] node[above] {$\hbox{ \sf\smalll stereographic projection }$}   
                                           node [below] {\sf\smalll {(star gazing)} }  (m-3-4)

        (m-1-1) edge  node[right]  {$\imath$}                              (m-2-1)
        (m-2-1) edge node[right] {$\pi $}
                             node[color=red, below =50pt]
                                     {\smalll\sf celestial sphere \ \ \ }     (m-3-1)
        (m-2-1) edge node[right] {$\pi $}                                     (m-3-1)
        (m-1-2) edge node[right] {$\subset$}                                (m-2-2)
         (m-1-3) edge [double]  node[right] {$\hbox{id}\otimes *$}             
                                      node[above=50pt, color=red] {$\hbox{\smalll\sf Pauli} \atop \hbox{\smalll\sf spinors}$}                      (m-2-3)
         (m-1-3.east) edge [dotted,bend left =145pt, out=60, in=150, very thick]  node[right] {$\pi$}                 (m-3-4)
         (m-1-3) edge [dashed, thin] node[above,  rotate=28] 
         {$\hbox{\smalll\sf parametrizing~~} \atop \hbox{\smalll\sf formula}$}      (m-2-2)
                             ;
         

\end{tikzpicture}   

{\small
~\\
{\bf Legend}\\[11pt]
$\dgaction{G}{~X~}$ \quad --- Group $G$ acting on set $X$, \ $X\ni x\mapsto gx\in X$\\[7pt]
$\dggaction{G}{~X~}$ \quad --- Group $G$ acting on set $X$, \ $X\ni x\mapsto gxg^*\in X$\\[5pt]
$\hbox{Herm}(2, \mathbb C)$ \quad --- two-by-two Hermitian matrices over $C$ \\[1pt]
$\hbox{Herm}_\circ(2, \mathbb C)$ \quad --- two-by-two traceless  Hermitian matrices over $C$  \\[1pt]
$\hbox{Herm}^\bullet(2, \mathbb C)$ \quad --- two-by-two   Hermitian matrices over $C$ of determinant equal to 0\\[1pt]
$\hbox{Symm}(2, \mathbb C)$ \quad --- two-by-two real symmetric matrices
 }

%
%
%
%
%
%
%
%
%


\section*{Acknowledgments}

The author would like to express his special thanks to Philip Feinsilver 
for very useful comments and encouragement.  
His effervescent attitude to mathematical explorations is priceless.


%
%
%
%
%

\end{document}